%% file: paper.tex
\documentclass[10pt, oneside, a4paper, abstracton]{scrartcl}

\usepackage[noadjust]{cite}

\usepackage{amsmath}
\usepackage{amssymb}
\usepackage{amsthm}

\usepackage{rotating}
\usepackage{multirow}

\usepackage{lmodern}
\usepackage{pifont}
\usepackage{dsfont}

\usepackage{tikz}
  \usetikzlibrary{calc}
  \usetikzlibrary{decorations.markings}
  \usetikzlibrary{decorations.pathmorphing}
  \usetikzlibrary{patterns}
  \usetikzlibrary{positioning}
  
\usepackage{pgfplots}
  \pgfplotsset{compat = 1.13}
  \usepgfplotslibrary{external}
  \tikzexternaldisable

\usepackage[linesnumbered, lined, algoruled]{algorithm2e}

\definecolor{myRed}{HTML}{E34A33}
\definecolor{myBlue}{HTML}{0571B0}
\definecolor{myBrown}{HTML}{A6611A}
\definecolor{linkBlue}{HTML}{0055C9}
\definecolor{linkRed}{HTML}{FF1A24}
\definecolor{linkPurple}{HTML}{6200D9}

\definecolor{plotcolor1}{HTML}{1B9E77}
\definecolor{plotcolor2}{HTML}{D95F02}
\definecolor{plotcolor3}{HTML}{7570B3}
\definecolor{plotcolor4}{HTML}{E7298A}
\definecolor{plotcolor5}{HTML}{66A61E}
\definecolor{plotcolor6}{HTML}{E6AB02}
\definecolor{plotcolor7}{HTML}{A6761D}
\definecolor{plotcolor8}{HTML}{666666}

\newcommand{\cA}{\ensuremath{\mathcal{A}}}
\newcommand{\cB}{\ensuremath{\mathcal{B}}}
\newcommand{\cC}{\ensuremath{\mathcal{C}}}
\newcommand{\cE}{\ensuremath{\mathcal{E}}}
\newcommand{\cH}{\ensuremath{\mathcal{H}}}
\newcommand{\cL}{\ensuremath{\mathcal{L}}}
\newcommand{\cQ}{\ensuremath{\mathcal{Q}}}
\newcommand{\cR}{\ensuremath{\mathcal{R}}}
\newcommand{\cT}{\ensuremath{\mathcal{T}}}
\newcommand{\cV}{\ensuremath{\mathcal{V}}}
\newcommand{\cU}{\ensuremath{\mathcal{U}}}
\newcommand{\hB}{\ensuremath{\widehat{B}}}
\newcommand{\hC}{\ensuremath{\widehat{C}}}
\newcommand{\hE}{\ensuremath{\widehat{E}}}
\newcommand{\hH}{\ensuremath{\widehat{H}}}
\newcommand{\hK}{\ensuremath{\widehat{K}}}
\newcommand{\hM}{\ensuremath{\widehat{M}}}
\newcommand{\hx}{\ensuremath{\hat{x}}}
\newcommand{\hy}{\ensuremath{\hat{y}}}
\newcommand{\hcA}{\ensuremath{\widehat{\mathcal{A}}}}
\newcommand{\hcB}{\ensuremath{\widehat{\mathcal{B}}}}
\newcommand{\hcC}{\ensuremath{\widehat{\mathcal{C}}}}
\newcommand{\hcE}{\ensuremath{\widehat{\mathcal{E}}}}
\newcommand{\tC}{\ensuremath{\widetilde{C}}}
\newcommand{\tE}{\ensuremath{\widetilde{E}}}
\newcommand{\tK}{\ensuremath{\widetilde{K}}}
\newcommand{\tP}{\ensuremath{\widetilde{P}}}
\newcommand{\tQ}{\ensuremath{\widetilde{Q}}}
\newcommand{\tq}{\ensuremath{\tilde{q}}}
\newcommand{\trans}{\ensuremath{\mkern-1.5mu\mathsf{T}}}
\newcommand{\herm}{\ensuremath{\mathsf{H}}}

\newcommand{\myP}{\textbf{\textit{p}}}
\newcommand{\myPM}{\textbf{\textit{pm}}}
\newcommand{\myPV}{\textbf{\textit{pv}}}
\newcommand{\myVP}{\textbf{\textit{vp}}}
\newcommand{\myVPM}{\textbf{\textit{vpm}}}
\newcommand{\myV}{\textbf{\textit{v}}}
\newcommand{\myFV}{\textbf{\textit{fv}}}
\newcommand{\mySO}{\textbf{\textit{so}}}

\DeclareMathOperator{\real}{Re}
\DeclareMathOperator{\imag}{Im}
\DeclareMathOperator{\pre}{pre}
\DeclareMathOperator{\diag}{diag}
\DeclareMathOperator{\range}{range}
\DeclareMathOperator{\mo}{mod}
\DeclareMathOperator{\step}{step}

\theoremstyle{definition}\newtheorem{definition}{Definition}
\theoremstyle{definition}\newtheorem{remark}{Remark}

\newcommand{\cmark}{\ding{51}}
\newcommand{\xmark}{\ding{55}}
\newcommand{\tabnum}[3]{{\small #1\texttt{e#2}#3}}

\pgfplotscreateplotcyclelist{freqlist}{
  {black, solid},
  {plotcolor1, mark options={solid}, mark=*, mark repeat=80, mark phase=0},
  {plotcolor2, mark options={solid}, mark=o, mark repeat=80, mark phase=8},
  {plotcolor3, mark options={solid}, mark=x, mark repeat=80, mark phase=16},
  {plotcolor4, mark options={solid}, mark=+, mark repeat=80, mark phase=24},
  {plotcolor5, mark options={solid}, mark=triangle*, mark repeat=80, mark phase=32},
  {plotcolor6, mark options={solid}, mark=square*, mark repeat=80, mark phase=40},
  {plotcolor7, mark options={solid}, mark=diamond*, mark repeat=80, mark phase=48},
  {plotcolor8, mark options={solid}, mark=Mercedes star, mark repeat=80, mark phase=56}
}

\pgfplotscreateplotcyclelist{timelist}{
  {black, solid},
  {plotcolor1, mark options={solid}, mark=*, mark repeat=400, mark phase=0},
  {plotcolor2, mark options={solid}, mark=o, mark repeat=400, mark phase=40},
  {plotcolor3, mark options={solid}, mark=x, mark repeat=400, mark phase=80},
  {plotcolor4, mark options={solid}, mark=+, mark repeat=400, mark phase=120},
  {plotcolor5, mark options={solid}, mark=triangle*, mark repeat=400, mark phase=160},
  {plotcolor6, mark options={solid}, mark=square*, mark repeat=400, mark phase=200},
  {plotcolor7, mark options={solid}, mark=diamond*, mark repeat=400, mark phase=240},
  {plotcolor8, mark options={solid}, mark=Mercedes star, mark repeat=400, mark phase=280}
}

\pgfplotsset{%
  compat=1.13,
  every axis plot/.append style={%
    line width = 1.0pt,
    mark size  = 2pt
  }
}

\tikzstyle{rngLine} = [black, line width = 1.5pt, dashed]

\usepackage[%
  colorlinks = true,
  linkcolor  = linkBlue,
  citecolor  = linkRed,
  urlcolor   = linkPurple
]{hyperref}
\usepackage{url}

\usepackage{geometry}
\usepackage{todonotes}

\begin{document}

\title{Frequency- and Time-Limited Balanced Truncation for Large-Scale
  Second-Order Systems}
\subtitle{Dedicated to Paul Van Dooren on the occasion of his 70th birthday.}
\author{%
  Peter~Benner\thanks{
   Max Planck Institute for Dynamics of Complex Technical Systems,
   Sandtorstr. 1, 39106 Magdeburg, Germany.\newline
   E-mail: \texttt{\href{mailto:benner@mpi-magdeburg.mpg.de}%
     {benner@mpi-magdeburg.mpg.de}}
   \newline
   Faculty of Mathematics, Otto von Guericke University,
   Universit{\"a}tsplatz 2, 39106 Magdeburg, Germany.\newline
   E-mail: \texttt{\href{mailto:peter.benner@ovgu.de}%
     {peter.benner@ovgu.de}}} \and
  Steffen~W.~R.~Werner\thanks{
    Max Planck Institute for Dynamics of Complex
    Technical Systems, Sandtorstr. 1, 39106 Magdeburg, Germany.\newline
    E-mail: \texttt{\href{mailto:werner@mpi-magdeburg.mpg.de}%
      {werner@mpi-magdeburg.mpg.de}}}
}
\date{~}

\maketitle


\begin{abstract}
  Considering the use of dynamical systems in practical applications, often
  only limited regions in the time or frequency domain are of interest.
  Therefor, it usually pays off to compute local approximations of the
  used dynamical systems in the frequency and time domain.
  In this paper, we consider a structure-preserving extension of the frequency-
  and time-limited balanced truncation methods to second-order dynamical
  systems.
  We give a full overview about the first-order limited balanced truncation
  methods and extend those methods to second-order systems by using the
  different second-order balanced truncation formulas from the literature.
  Also, we present numerical methods for solving the arising large-scale sparse
  matrix equations and give numerical modifications to deal with the problematic
  case of second-order systems.
  The results are then illustrated on three numerical examples.
 
  \vspace{1em}
  \textbf{Keywords:} model order reduction, second-order differential equations,
    linear systems, balanced truncation, frequency-limited balanced truncation,
    time-limited balanced truncation, local model reduction,
    structure-preserving approximation
\end{abstract}


\section{Introduction}

The modeling of, e.g., mechanical and electrical systems often leads to
linear dynamical systems containing second-order time derivatives.
In this paper, we consider linear second-order input-output systems of
the form
\begin{align} \label{eqn:sosys}
  \begin{aligned}
    M\ddot{x} + E\dot{x}(t) + Kx(t) & = B_{u}u(t),\\
    y(t) & = C_{p}x(t) + C_{v}\dot{x}(t),
  \end{aligned}
\end{align}
with $M, E, K \in \mathbb{R}^{n \times n}$, $B_{u} \in \mathbb{R}^{n \times m}$
and $C_{p}, C_{v} \in \mathbb{R}^{p \times n}$, and $u(t) \in \mathbb{R}^{m}$,
the inputs, $x(t) \in \mathbb{R}^{n}$, the states, and
$y(t) \in \mathbb{R}^{p}$, the outputs of the system.
In the frequency domain, the input-to-output relation is directly given
as $y(s) = H(s)u(s)$, whereby the so-called transfer function is given by
\begin{align} \label{eqn:sotf}
  H(s) & = (sC_{v} + C_{p})(s^{2}M + sE + K)^{-1}B_{u},
\end{align}
with $s \in \mathbb{C}$.
In applications, the number of differential equations, $n$, describing the
system, can become very large.
This complicates using the model for simulations and controller design due to
the expensive costs in terms of computational resources as time and memory.
Therefor, model reduction is needed to construct a surrogate system with
a much smaller number of equations $r \ll n$, which approximates the
input-to-output behavior of~\eqref{eqn:sosys}.
To use the surrogate model as the original one, e.g., applying the same tools,
the surrogate needs to have the same structure as the original system, i.e.,
the reduced-order model should also have the form
\begin{align*}
  \begin{aligned}
    \hM \ddot{\hx}(t) + \hE \dot{\hx}(t) + \hK \hx(t) & = \hB_{u}u(t),\\
    \hy(t) & = \hC_{p}\hx(t) + \hC_{v}\dot{\hx}(t),
  \end{aligned}
\end{align*}
with the new system matrices $\hM, \hE, \hK \in \mathbb{R}^{r \times r}$,
$\hB_{u} \in \mathbb{R}^{r \times m}$ and $\hC_{p}, \hC_{v} \in
\mathbb{R}^{p \times r}$.

Due to its relevance in a lot of applications, the problem of
structure-preserving model reduction for second-order systems
has already been investigated in the literature to quiet an extend.
There are structure-preserving extensions of classical model reduction methods
like modal truncation and dominant pole algorithms~\cite{RomM06, morBenKTetal16,
morSaaSW19}, moment matching~\cite{morBeaG05, morSal05, morChaGVetal05,
morSalL06}, balanced truncation~\cite{morMeyS96, morChaLVetal06, morReiS08}, or
for example of the $\mathcal{H}_{2}$-optimal iterative rational Krylov
algorithm~\cite{morWya12}.
Especially, we want to mention the work of Paul Van Dooren, and co-authors, on
the second-order balanced truncation approach.
In~\cite{morChaGVetal05}, he introduced a new balancing idea that is stronger
related to the origins of balanced truncation than the other extensions.
Most of the extended methods aim for a globally good approximation behavior,
but very often, only the local system's behavior in the frequency or time domain
is of actual interest for the application.
In case of first-order systems, the frequency- and time-limited balanced
truncation methods, first mentioned in~\cite{morGawJ90}, aiming for such local
approximations.
Those methods have been extended in the first-order case to large-scale
sparse systems~\cite{morBenKS16, morKue18} and to system with
differential-algebraic equations~\cite{morImrG15, morHaiGIetal17}.

A first attempt to generalize the limited balanced truncation methods to
second-order systems has been done in~\cite{morHaiGIetal18} for the
frequency-limited balanced truncation by making use of some formulas
from~\cite{morReiS08} and for the time-limited balanced truncation
in~\cite{morHaiGIetal18a} in the same way.
In this paper, we are extending the frequency- and time-limited balanced
truncation methods by using all the different second-order balanced truncation
approaches from the literature~\cite{morMeyS96, morReiS08, morChaLVetal06}
and correct some mistakes that were made
in~\cite{morHaiGIetal18, morHaiGIetal18a} considering the issue of stability
preservation.
Also, we are extending the numerical approaches to the large-scale second-order
system case and present strategies to deal with numerical difficulties
aligning with second-order systems in general.

The paper has the following structure.
Section~\ref{sec:lbt} contains a review of the theory for the classical and
limited balanced truncation methods in the generalized first-order system
case; see Section~\ref{sec:fo}; as well as a review of the different
second-order balanced approaches and the extensions of the limited balanced
truncation methods to second-order systems in Section~\ref{sec:so}.
Afterwards, in Section~\ref{sec:nummeth}, the numerical methods for solving the
large-scale sparse matrix equations with function right hand-side are covered.
Also, the $\alpha$-shift strategy and two-step methods are explained in this
section, which ends with the modified Gramian approach and remarks on the
stability preservation of the methods.
Three numerical examples are then given in Section~\ref{sec:examples} to
demonstrate the applicability of the methods on large-scale sparse second-order
systems. Section~\ref{sec:conclusions} concludes the paper.


\section{The frequency- and time-limited balanced truncation methods}
\label{sec:lbt}


\subsection{First-order system case}%
\label{sec:fo}

In this section, we will remind of the classical balanced truncation technique
and give an overview on the frequency- and time-limited versions of this method
for the case of first-order systems.


\subsubsection{Classical balanced truncation}%
\label{sec:fobt}

We consider here generalized first-order state-space systems of the form
\begin{align} \label{eqn:fosys}
  \begin{aligned}
    \cE\dot{q}(t) & = \cA q(t) + \cB u(t),\\
    y(t) & = \cC q(t),
  \end{aligned}
\end{align}
with $\cE, \cA \in \mathbb{R}^{N \times N}$, $\cB \in \mathbb{R}^{N \times m}$,
$\cC \in \mathbb{R}^{p \times N}$, and the corresponding transfer function
\begin{align} \label{eqn:fotf}
  H(s) & = \cC(s\cE - \cA)^{-1}\cB,
\end{align}
with $s \in \mathbb{C}$.
For simplicity, we are assuming that $\mathcal{E}$ is invertible and the system
is c-stable, i.e., all eigenvalues of $\lambda\cE - \cA$ lie in the open left
complex half-plane.
The extension of the balanced truncation method to the descriptor system case
($\cE$ non-invertible) can be found in~\cite{Sty02, morBenS17}.
The system Gramians of~\eqref{eqn:fosys} are
defined as
\begin{align}\label{eqn:fogram}
  \begin{aligned}
  P_{\infty} & = \frac{1}{2\pi}\int\limits_{-\infty}^{+\infty}
    {(j\omega \cE - \cA)^{-1}\cB\cB^{\trans}(-j\omega \cE - \cA)^{-\trans}
    \mathrm{d}\omega}
    = \int\limits_{0}^{+\infty}{e^{\cE^{-1}\cA t}\cE^{-1}\cB\cB^{\trans}
    \cE^{-\trans}e^{\cA^{\trans}\cE^{-\trans} t}\mathrm{d}t},\\
  Q_{\infty} & = \frac{1}{2\pi}\int\limits_{-\infty}^{+\infty}
    {(-j\omega \cE - \cA)^{-\trans}\cC^{\trans}\cC(j\omega \cE - \cA)^{-1}
    \mathrm{d}\omega} = \int\limits_{0}^{+\infty}{\cE^{-\trans}
    e^{\cA^{\trans}\cE^{-\trans} t}\cC^{\trans}\cC e^{\cE^{-1}\cA t}\cE^{-1}
    \mathrm{d}t},
  \end{aligned}
\end{align}
with $P_{\infty}$ the infinite controllability Gramian and
$\cE^{\trans} Q_{\infty}\cE$ the infinite observability Gramian.
Note that in the infinite case, the frequency and time representations of the 
Gramians are equal.
It can be shown that those Gramians~\eqref{eqn:fogram} are the unique, 
symmetric positive semi-definite solutions of the following Lyapunov equations
\begin{align} \label{eqn:folyap}
  \begin{aligned}
  \cA P_{\infty} \cE^{\trans} + \cE P_{\infty}\cA^{\trans}
    + \cB \cB^{\trans} & = 0,\\
  \cA^{\trans} Q_{\infty} \cE + \cE^{\trans} Q_{\infty} \cA
    + \cC^{\trans} \cC & = 0.
  \end{aligned}
\end{align}
The Hankel singular values are then defined as the positive square-roots
of the eigenvalues of $P_{\infty}\cE^{\trans}Q_{\infty}\cE$, which
are a measure of how much influence the corresponding states have on the
input-output behavior of the system.
The main idea of balanced truncation is to balance the system such that
\begin{align*}
  P_{\infty} & = Q_{\infty} = \begin{bmatrix} \sigma_{1} & & & \\
    & \sigma_{2} & & \\ & & \ddots & \\ & & & \sigma_{N} \end{bmatrix},
\end{align*}
with the Hankel singular values $\sigma_{1} \geq \sigma_{2} \geq \ldots \geq
\sigma_{N} > 0$ and then to truncate states corresponding to small Hankel
singular values~\cite{morMoo81}.
The complete balanced truncation square-root method is summarized in
Algorithm~\ref{alg:bt}.

\begin{algorithm}[t]
  \DontPrintSemicolon
  \caption{Balanced Truncation Square-Root Method}
  \label{alg:bt}
  
  \KwIn{System matrices $\cA$, $\cB$, $\cC$, $\cE$ from~\eqref{eqn:fosys}.}
  \KwOut{Matrices of the reduced-order system $\hcA$, $\hcB$, $\hcC$,
    $\hcE$.}
  
  Compute Cholesky factorizations of the Gramians by solving the Lyapunov
    equations~\eqref{eqn:folyap} such that
    $P_{\infty} = R_{\infty}R_{\infty}^{\trans}$,
    $Q_{\infty} = L_{\infty}L_{\infty}^{\trans}$.\;
  Compute the singular value decomposition
    \begin{align*}
      L_{\infty}^{\trans} \cE R_{\infty} & =
        \begin{bmatrix} U_{1} & U_{2} \end{bmatrix}
        \begin{bmatrix} \Sigma_{1} & \\ & \Sigma_{1} \end{bmatrix}
        \begin{bmatrix} V_{1}^{\trans} \\ V_{2}^{\trans} \end{bmatrix},
    \end{align*}
    with $\Sigma_{1} = \mathrm{diag}(\sigma_{1}, \ldots, \sigma_{r})$ containing
    the $r$ largest Hankel singular values.\;
  Construct the projection matrices
    \begin{align*}
      \begin{aligned}
        T & = R_{\infty}V_{1}\Sigma_{1}^{-\frac{1}{2}} && \text{and} &
          W & = L_{\infty}U_{1}\Sigma_{1}^{-\frac{1}{2}}.
      \end{aligned}
    \end{align*}\vspace{-\baselineskip}\;
  Compute the reduced-order model by
    \begin{align*}
      \begin{aligned}
        \hcA & = W^{\trans}\cA T, & \hcB & = W^{\trans}\cB, &
          \hcC & = \cC T, & \hcE & = W^{\trans}\cE T = I_{r}.
      \end{aligned}
    \end{align*}\vspace{-\baselineskip}\;
\end{algorithm}

The balanced truncation method provides an a posteriori error bound in the
$\cH_{\infty}$ norm
\begin{align} \label{eqn:btbound}
  \lVert H - \hH \rVert_{\cH_{\infty}} \leq 2\sum\limits_{k = r + 1}^{N}
    {\sigma_{k}^{2}},
\end{align}
where $H$ is the transfer function of the original model~\eqref{eqn:fotf} and
$\hH$ the transfer function of the reduced-order model.
The bound~\eqref{eqn:btbound} depends only on the truncated Hankel singular
values, which allows an adaptive choice of the reduction order.
Also, this method preserves the stability of the original model, i.e., if $H$
was a c-stable model then also $\hH$ will be c-stable.

The application of the balanced truncation method to large-scale sparse systems
is possible by approximating the Cholesky factors of the Gramians via low-rank
factors $P_{\infty} \approx Z_{R_{\infty}}Z_{R_{\infty}}^{\trans}$,
$Q_{\infty} \approx Z_{L_{\infty}}Z_{L_{\infty}}^{\trans}$, with
$Z_{R_{\infty}} \in \mathbb{R}^{N \times k_{R}}$,
$Z_{L_{\infty}} \in \mathbb{R}^{N \times k_{L}}$ and
$k_{R}, k_{L} \ll N$; see, e.g.,~\cite{morBenQQ00}.
The approximation of the Gramians is reasonable due to a fast singular value
decay arising by the low-rank right-hand sides~\cite{BakES15}.
For the computation of those factors, appropriate low-rank techniques are well
developed~\cite{BenS13}.


\subsubsection{Frequency-limited approach}%
\label{sec:foflbt}

A suitable method to localize the approximation behavior of the balanced
truncation method in the frequency domain is the frequency-limited balanced
truncation~\cite{morGawJ90}.
The idea is based on the frequency representation of the system 
Gramians~\eqref{eqn:fogram}, such that the frequency-limited Gramians 
of~\eqref{eqn:fosys} are given by
\begin{align}\label{eqn:foflgram}
  \begin{aligned}
  P_{\Omega} = \frac{1}{2\pi} \int\limits_{\Omega}
    {(j\omega \cE - \cA)^{-1}\cB\cB^{\trans}(-j\omega \cE - \cA)^{-\trans}
    \mathrm{d}\omega},\\
  Q_{\Omega} = \frac{1}{2\pi} \int\limits_{\Omega}
    {(-j\omega \cE - \cA)^{-\trans}\cC^{\trans}\cC(j\omega \cE - \cA)^{-1}
    \mathrm{d}\omega},
  \end{aligned}
\end{align}
where $\Omega = [-\omega_{2}, -\omega_{1}] \cup [\omega_{1}, \omega_{2}] \subset 
\mathbb{R} $ is the frequency range of interest.
It can be shown that the left-hand sides in~\eqref{eqn:foflgram} are also given
as the unique, symmetric positive semi-definite solutions of the two Lyapunov
equations
\begin{align}\label{eqn:fofllyap}
  \begin{aligned}
    \cA P_{\Omega} \cE^{\trans} + \cE P_{\Omega} \cA^{\trans}
      + \cB_{\Omega}\cB^{\trans} + \cB\cB_{\Omega}^{\trans} & = 0,\\
    \cA^{\trans} Q_{\Omega} \cE + \cE^{\trans} Q_{\Omega} \cA
      + \cC_{\Omega}^{\trans}\cC + \cC^{\trans}\cC_{\Omega} & = 0,
  \end{aligned}
\end{align}
with new right hand-side matrices $B_{\Omega} = \cE F_{\Omega} \cB$,
$C_{\Omega} = \cC F_{\Omega} \cE$ containing the matrix functions
\begin{align} \label{eqn:fomega}
  \begin{aligned}
    F_{\Omega} & = \real\left(\frac{j}{\pi}\ln\left((\cA + j\omega_{1}\cE)^{-1} 
      (\cA + j\omega_{2}\cE)\right)\right)\cE^{-1}\\
    & = \cE^{-1} \real\left(\frac{j}{\pi}\ln\left((\cA + j\omega_{2}\cE)
      (\cA + j\omega_{1}\cE)^{-1}\right)\right),
  \end{aligned}
\end{align}
with $\ln(.)$ the principle branch of the matrix logarithm. 
Note that in case of $\Omega = [-\omega, \omega]$, the function 
evaluation~\eqref{eqn:fomega} simplifies to
\begin{align*}
  F_{\Omega} & = \real\left(\frac{j}{\pi}\ln\left(-\cE^{-1}\cA - j\omega I_{n} 
    \right)\right)\cE^{-1}\\
  & = \cE^{-1}\real\left(\frac{j}{\pi}\ln\left(-\cA\cE^{-1} - j\omega I_{n} 
    \right)\right).
\end{align*}
Also, the frequency-limited Gramians can be extended to an arbitrary number of
frequency intervals, i.e., for 
\begin{align*}
  \Omega = \bigcup\limits_{k = 1}^{\ell}\left([-\omega_{2k}, \omega_{2k-1}] 
    \cup [\omega_{2k-1}, \omega_{2k}]\right),
\end{align*}
with $0 < \omega_{1} < \ldots < \omega_{\ell}$, leads to the following 
modification of~\eqref{eqn:fomega}
\begin{align*}
  F_{\Omega} & = \real \left( \frac{j}{\pi} \ln \left(
    \prod\limits_{k = 1}^{\ell} (A + j\omega_{2k-1}\cE)^{-1}
    (A + j\omega_{2k}\cE) \right)\right)\cE^{-1}\\
  & = \cE^{-1}\real \left( \frac{j}{\pi} \ln \left(
      \prod\limits_{k = 1}^{\ell} (A + j\omega_{2k}\cE)
      (A + j\omega_{2k-1}\cE)^{-1}\right)\right).
\end{align*}
See~\cite{morBenKS16} for a more detailed discussion of the theory addressed
above.
The extension of this method to the large-scale system case can also be found
in~\cite{morBenKS16} and an extension to descriptor systems in~\cite{morImrG15}.
The resulting frequency-limited balanced truncation method is summarized in 
Algorithm~\ref{alg:flbt}.

\begin{algorithm}[t]
  \DontPrintSemicolon
  \caption{Frequency-Limited Balanced Truncation Square-Root Method}
  \label{alg:flbt}
  
  \KwIn{System matrices $\cA$, $\cB$, $\cC$, $\cE$ from~\eqref{eqn:fosys},
    frequency range of interest $\Omega$.}
  \KwOut{Matrices of the reduced-order system $\hcA$, $\hcB$, $\hcC$,
    $\hcE$.}
    
  Compute Cholesky factorizations of the frequency-limited Gramians 
    by solving the frequency-limited Lyapunov equations~\eqref{eqn:fofllyap}
    such that $P_{\Omega} = R_{\Omega}R_{\Omega}^{\trans}$,
    $Q_{\Omega} = L_{\Omega}L_{\Omega}^{\trans}$.\;
    
  Follow the steps 2--4 in Algorithm~\ref{alg:bt}.\;
\end{algorithm}


\subsubsection{Time-limited approach}%
\label{sec:fotlbt}

The counterpart of the frequency-limited balanced truncation from the previous
section in the time domain is the time-limited balanced
truncation~\cite{morGawJ90}.
This method aims for the approximation of the system on a time interval $T =
[t_{0}, t_{f}]$, where $0 \leq t_{0} < t_{f}$, based on the limitation of the 
time domain representation of the Gramians~\eqref{eqn:fogram}.
The time-limited Gramians of~\eqref{eqn:fosys} are then given by
\begin{align}\label{eqn:fotlgram}
  \begin{aligned}
    P_{T} & = \int\limits_{t_{0}}^{t_{f}}
      e^{\cE^{-1}\cA t}\cE^{-1}\cB\cB^{\trans}\cE^{-\trans}
      e^{\cA^{\trans}\cE^{-\trans} t}\mathrm{d}t,\\
    Q_{T} & = \int\limits_{t_{0}}^{t_{f}}
      \cE^{-\trans}e^{\cA^{\trans}\cE^{-\trans} t}\cC^{\trans}
      \cC e^{\cE^{-1}\cA t}\cE^{-1}\mathrm{d}t.
  \end{aligned}
\end{align}
and it can be shown, that the left-hand sides in~\eqref{eqn:fotlgram} are the
unique, positive semi-definite solutions of the two following Lyapunov equations
\begin{align}\label{eqn:fotllyap}
  \begin{aligned}
    \cA P_{T} \cE^{\trans} + \cE P_{T} \cA^{\trans}
      + \cB_{t_{0}}\cB_{t_{0}}^{\trans}
      - \cB_{t_{f}}\cB_{t_{f}}^{\trans} & = 0,\\
    \cA^{\trans} Q_{T} \cE + \cE^{\trans} Q_{T} \cA
      + \cC_{t_{0}}^{\trans}\cC_{t_{0}}
      - \cC_{t_{f}}^{\trans}\cC_{t_{f}} & = 0,
  \end{aligned}
\end{align}
where the new right hand-side matrices $\cB_{t_{0/f}} = \cE
e^{\cE^{-1}At_{0/f}}\cE^{-1}\cB = e^{A\cE^{-1}t_{0/f}}\cB$ and
$\cC_{t_{0/f}} = \cC e^{\cE^{-1}At_{0/f}}$ contain the matrix exponential.
The right hand-sides of~\eqref{eqn:fotllyap} simplify in case of
$t_{0} = 0$ since $\cB_{0} = \cB$ and $\cC_{0} = \cC$.
A more detailed discussion of the time-limited theory, especially for the 
large-scale system case, can be found in~\cite{morKue18}.
Also, the extension of the theory to the case of descriptor systems is given
in~\cite{morHaiGIetal17}.
It can be noted that considering more than one time interval at once
$[t_{0,1},t_{f,1}] \cup \cdots \cup [t_{0,\ell},t_{f,\ell}]$ is not
practical and usually one cannot guarantee a good approximation behavior in the
single intervals.
Instead it is common to take the smallest and largest time points in the 
intervals to construct a new overarching time interval
$[t_{0,\min}, t_{f,\max}]$, where $t_{0,\min} = \min\{t_{0,1}, \ldots,
t_{0,\ell}\}$ and $t_{0,\max} = \max\{t_{f,1}, \ldots, t_{f,\ell}\}$ such that
\begin{align*}
  \bigcup\limits_{k = 1}^{\ell}[t_{0,k}, t_{f,k}] \subset
    [t_{0,\min}, t_{f,\max}] = T.
\end{align*}
The resulting time-limited balanced truncation method is summarized in
Algorithm~\ref{alg:tlbt}.

\begin{algorithm}[t]
  \DontPrintSemicolon
  \caption{Time-Limited Balanced Truncation Square-Root Method}
  \label{alg:tlbt}
  
  \KwIn{System matrices $\cA$, $\cB$, $\cC$, $\cE$ from~\eqref{eqn:fosys},
    time range of interest $T$.}
  \KwOut{Matrices of the reduced-order system $\hcA$, $\hcB$, $\hcC$,
    $\hcE$.}
    
  Compute Cholesky factorizations of the time-limited Gramians 
    by solving the time-limited Lyapunov equations~\eqref{eqn:fotllyap} such that
    $P_{T} = R_{T}R_{T}^{\trans}$,
    $Q_{T} = L_{T}L_{T}^{\trans}$.\;
    
  Follow the steps 2--4 in Algorithm~\ref{alg:bt}.\;
\end{algorithm}


\subsection{Second-order case}%
\label{sec:so}

After recapitulating the basic ideas of the classical as well as the frequency-
and time-limited balanced truncation methods for first-order systems,
in this section we will extend those methods to second-order 
systems~\eqref{eqn:sosys}.


\subsubsection{Second-order balanced truncation methods}%
\label{sec:sobt}

Over time, there have been many attempts for the generalization of the 
classical balanced truncation method to the second-order
system case~\cite{morReiS08, morMeyS96, morChaLVetal06}.
All of them have in common the idea of linearization, i.e., the second-order
system~\eqref{eqn:sosys} is rewritten as a first-order system.
The usual linearization of choice for~\eqref{eqn:sosys} is its so-called first
companion form
\begin{align} \label{eqn:socomp}
  \begin{aligned}
  \underbrace{\begin{bmatrix} J & 0 \\ 0 & M \end{bmatrix}}_{\cE} \dot{q}(t) & =
    \underbrace{\begin{bmatrix} 0 & J \\ -K & -E \end{bmatrix}}_{\cA} q(t)
    + \underbrace{\begin{bmatrix} 0 \\ B_{u} \end{bmatrix}}_{\cB},\\
  y(t) & = \underbrace{\begin{bmatrix} C_{p} & C_{v} \end{bmatrix}}_{\cC}q(t),
  \end{aligned}
\end{align}
where $q(t) = [x^{\trans}(t), \dot{x}^{\trans}(t)]^{\trans}$ is the new combined 
state vector. 
The matrix $J \in \mathbb{R}^{n \times n}$ is an arbitrary invertible matrix but
usually chosen as $J = I_{n}$ or $J = -K$, which can lead to symmetric $\cA$ and
$\cE$ matrices in case of mechanical systems.

For system~\eqref{eqn:socomp}, the first-order Gramians are used, as given
by~\eqref{eqn:fogram} or~\eqref{eqn:folyap}, and then partitioned
according to the block structure in~\eqref{eqn:socomp} such that
\begin{align}\label{eqn:sogrampart}
  \begin{aligned}
    P_{\infty} & = \begin{bmatrix} P_{p} & P_{12} \\
      P_{12}^{\trans} & P_{v} \end{bmatrix}
      & \text{and} && Q_{\infty} & = \begin{bmatrix} Q_{p} & Q_{12} \\ 
      Q_{12}^{\trans} & Q_{v} \end{bmatrix},
  \end{aligned}
\end{align}
where $P_{p}$, $Q_{p}$ are the the position Gramians 
of~\eqref{eqn:sosys} and $P_{v}$, $Q_{v}$ the velocity Gramians.
Due to $P_{\infty} = P_{\infty}^{\trans} \geq 0$ and
$Q_{\infty} = Q_{\infty}^{\trans} \geq 0$, also the position and
velocity Gramians are symmetric positive semi-definite and can be written
in terms of their Cholesky factorizations
\begin{align*}
  \begin{aligned}
    P_{p} & = R_{p}R_{p}^{\trans}, & P_{v} & = R_{v}R_{v}^{\trans},
      & Q_{p} & = L_{p}L_{p}^{\trans}, & Q_{v} & = L_{v}L_{v}^{\trans}.
  \end{aligned}
\end{align*}
Based on those, the different second-order balanced truncation methods are
defined by balancing certain combinations of the four position and velocity
Gramians.
For most of the methods, the resulting balanced truncation is computed as
second-order projection method
\begin{align} \label{eqn:soproj}
  \begin{aligned}
    \hM & = WMT, & \hE & = WET, & \hK & = WKT, & \hB_{u} & = WB_{u}, & \hC_{p}
      & = C_{p}T, & \hC_{v} & = C_{v}T,
  \end{aligned}
\end{align}
where the different choices for $W$ and $T$ can be found in
Table~\ref{tab:sobt}.
There, the different transformation formulas are summarized and denoted
by the type as used in the corresponding references.
The subscript $1$ matrices denote the part of the singular value
decompositions corresponding to the largest characteristic singular values.

\begin{table}
  \centering
  \begin{tabular}{c|l|l|c}
    \multicolumn{1}{c|}{\textbf{Type}}
      & \multicolumn{1}{c|}{\textbf{SVD(s)}}
      & \multicolumn{1}{c|}{\textbf{Transformation}}
      & \multicolumn{1}{c}{\textbf{Reference}}\\ \hline
    \textbf{\myV}
      & $U \Sigma V^{\trans} = L_{v}^{\trans} M R_{v}$
      & $W = L_{v}U_{1}\Sigma_{1}^{-\frac{1}{2}},~
        T = R_{v}V_{1}\Sigma_{1}^{-\frac{1}{2}}$
      & \cite{morReiS08}%
        \vphantom{\Large $P^{\frac{1}{2}}_{\frac{1}{2}}$}\\ \hline
    \textbf{\myFV}
      & $\ast\, \Sigma V^{\trans} = L_{p}^{\trans} J R_{p}$
      & $W = T,~
         T = R_{p}V_{1}\Sigma_{1}^{-\frac{1}{2}}$
      & \cite{morMeyS96}%
        \vphantom{\Large $P^{\frac{1}{2}}_{\frac{1}{2}}$}\\ \hline
    \textbf{\myVPM}
      & $U \Sigma V^{\trans} = L_{p}^{\trans} J R_{v}$
      & $W = M^{-\trans} J^{\trans} L_{p}U_{1}\Sigma_{1}^{-\frac{1}{2}},~
         T = R_{v}V_{1}\Sigma_{1}^{-\frac{1}{2}}$
      & \cite{morReiS08}%
        \vphantom{\Large $P^{\frac{1}{2}}_{\frac{1}{2}}$}\\ \hline
    \textbf{\myPM}
      & {\arraycolsep = 1.4pt
        $\begin{array}{rcl}
           U \Sigma V^{\trans} & = & L_{p}^{\trans} J R_{p},
             \vphantom{\text{\Large} P^{\frac{1}{2}}_{1}}\\
         \end{array}$}
      & $W = M^{-\trans} J^{\trans} L_{p}U_{1}\Sigma_{1}^{-\frac{1}{2}},~
         T = R_{p}V_{1}\Sigma_{1}^{-\frac{1}{2}}$
      & \cite{morReiS08}%
        \vphantom{\Large $P^{\frac{1}{2}}_{\frac{1}{2}}$}\\ \hline
    \textbf{\myPV}
      & $U \Sigma V^{\trans} = L_{v}^{\trans} M R_{p}$
      & $W = L_{v}U_{1}\Sigma_{1}^{-\frac{1}{2}},~
         T = R_{p}V_{1}\Sigma_{1}^{-\frac{1}{2}}$
      & \cite{morReiS08}%
        \vphantom{\Large $P^{\frac{1}{2}}_{\frac{1}{2}}$}\\ \hline
    \textbf{\myVP}
      & {\arraycolsep = 1.4pt
        $\begin{array}{rcl}
           \ast \Sigma V^{\trans} & = & L_{p}^{\trans} J R_{v},
             \vphantom{\text{\Large} P^{\frac{1}{2}}_{1}}\\
           U \ast \ast & = & L_{v}^{\trans} M R_{p}
             \vphantom{\text{\Large} P^{\frac{1}{2}}_{\frac{1}{2}}}
         \end{array}$}
      & $W = L_{v}U_{1}\Sigma_{1}^{-\frac{1}{2}},~
         T = R_{v}V_{1}\Sigma_{1}^{-\frac{1}{2}}$
      & \cite{morReiS08}\\ \hline
    \textbf{\myP}
      & {\arraycolsep = 1.4pt
        $\begin{array}{rcl}
           \ast\, \Sigma V^{\trans} & = & L_{p}^{\trans} J R_{p},
             \vphantom{\text{\Large} P^{\frac{1}{2}}_{1}}\\
           U \ast \ast & = & L_{v}^{\trans} M R_{v}
             \vphantom{\text{\Large} P^{\frac{1}{2}}_{\frac{1}{2}}}
         \end{array}$}
      & $W = L_{v}U_{1}\Sigma_{1}^{-\frac{1}{2}},~
         T = R_{p}V_{1}\Sigma_{1}^{-\frac{1}{2}}$
      & \cite{morReiS08}\\ \hline
    \textbf{\mySO}
      & {\arraycolsep = 1.4pt
        $\begin{array}{rcl}
          U_{p} \Sigma_{p} V_{p}^{\trans} & = & L_{p}^{\trans} J R_{p},
            \vphantom{\text{\Large} P^{\frac{1}{2}}_{1}}\\
          U_{v} \Sigma_{v} V_{v} & = & L_{v}^{\trans} M R_{v}
            \vphantom{\text{\Large} P^{\frac{1}{2}}_{\frac{1}{2}}}
         \end{array}$}
      & {\arraycolsep = 1.4pt
        $\begin{array}{rclrcl}
           W_{p} & = & L_{p}U_{p,1}\Sigma_{p,1}^{-\frac{1}{2}},~&
           T_{p} & = & R_{p}V_{p,1}\Sigma_{p,1}^{-\frac{1}{2}},\\
           W_{v} & = & L_{v}U_{v,1}\Sigma_{v,1}^{-\frac{1}{2}},~&
           T_{v} & = & R_{v}V_{v,1}\Sigma_{v,1}^{-\frac{1}{2}}
         \end{array}$}
      & \cite{morChaLVetal06}\\ \hline
  \end{tabular}
  \caption{Second-order balanced truncation formulas. (Here, $\ast$ denotes
    factors of the SVD not needed, and thus not accumulated in practical
    computations.)}
  \label{tab:sobt}
\end{table}

In contrast to the balancing methods that describe the reduced-order model
by~\eqref{eqn:soproj}, the second-order balanced truncation (\mySO)
from~\cite{morChaLVetal06} computes the reduced-order model by
\begin{align}\label{eqn:sobtproj}
  \begin{aligned}
    \hM  & = S \left( W_{v}^{\trans} M T_{v} \right) S^{-1}, &
      \hE & = S \left( W_{v}^{\trans} E T_{v} \right) S^{-1}, &
      \hK & = S \left( W_{v}^{\trans} K T_{p} \right),\\
    \hB_{u} & = S \left( W_{v}^{\trans} B_{u} \right), &
      \hC_{p} & = C_{p} T_{p}, &
      \hC_{v} & = C_{v} T_{v} S^{-1},
  \end{aligned}
\end{align}
where $S = W_{p} J T_{v}$ and the transformation matrices $W_{p}$, $W_{v}$,
$T_{p}$, $T_{v}$ are given in the last line of Table~\ref{tab:sobt}.
This type of balancing can be seen as a projection method for the first-order
realization~\eqref{eqn:socomp} with a recovering of the second-order structure.

The general second-order balanced truncation square-root method is summarized
in Algorithm~\ref{alg:sobt}.

\begin{remark}
  In contrast to the first-order balanced truncation described in 
  Section~\ref{sec:fobt}, none of the second-order balanced truncation
  methods provides an error bound in the $\mathcal{H}_{\infty}$ norm
  or can preserve the stability of the original system in the general case.
  A collection of examples for the stability issue is given in~\cite{morReiS08}.
  In case of symmetric second-order systems, i.e., $M = M^{\trans}$,
  $E = E^{\trans}$, $K = K^{\trans}$, $C_{p} = B_{u}^{\trans}$ and $C_{v} = 0$,
  it can be shown that the position-velocity balancing (\myPV) as well as the
  free-velocity balancing (\myFV) are stability preserving.
  Note that the position-velocity balancing also belongs to the class of
  balanced truncation approaches, which define system Gramians by using the
  underlying transfer function structure~\eqref{eqn:sotf}.
  Those balancing approaches have been generalized in~\cite{morBre16} for 
  systems with integro-differential equations.
\end{remark}

\begin{algorithm}[t]
  \DontPrintSemicolon
  \caption{Second-Order Balanced Truncation Square-Root Method}
  \label{alg:sobt}
  
  \KwIn{System matrices $M$, $E$, $K$, $B_{u}$, $C_{p}$, $C_{v}$
    from~\eqref{eqn:sosys}.}
  \KwOut{Matrices of the reduced-order system $\hM$, $\hE$, $\hK$, $\hB_{u}$
    $\hC_{p}$, $\hC_{v}$.}
    
  Compute Cholesky factorizations of the first-order system Gramians 
    by solving~\eqref{eqn:folyap}, where the linearization~\eqref{eqn:socomp}
    is used, such that
    $P_{\infty} = R_{\infty}R_{\infty}^{\trans}$,
    $Q_{\infty} = L_{\infty}L_{\infty}^{\trans}$.\;
    
  Partition the Cholesky factors according to the first-order formulation
    \begin{align*}
      \begin{aligned}
        R_{\infty} & = \begin{bmatrix} R_{p} \\ R_{v} \end{bmatrix} &
          \text{and} && L_{\infty} & = \begin{bmatrix} L_{p} \\
          L_{v} \end{bmatrix}.
      \end{aligned}
    \end{align*}
    \vspace{-\baselineskip}\;
    
  Compute the singular value decompositions and transformation matrices
    as in Table~\ref{tab:sobt}.\;
    
  Compute the reduced-order model by either~\eqref{eqn:soproj} for the methods
    \myP, \myPM, \myPV, \myVP, \myVPM, \myV{} and \myFV{} or
    by~\eqref{eqn:sobtproj} for \mySO.\;
\end{algorithm}


\subsubsection{Second-order frequency-limited approach}%
\label{sec:soflbt}

The generalization of the frequency-limited balanced truncation method for
second-order systems has been discussed in~\cite{morHaiGIetal18} for the
position (\myP) and position-velocity (\myPV) balancing from~\cite{morReiS08}.
Here we will summarize their results and give a more general extension for
the frequency-limited second-order balanced truncation method.
The basic idea for the approach comes from the observation that the block
partitioning of the Gramians~\eqref{eqn:sogrampart} can be written as
\begin{align}\label{eqn:sogramproj}
  \begin{aligned}
    P_{p} & = \begin{bmatrix} I_{n} & 0 \end{bmatrix} P_{\infty}
      \begin{bmatrix} I_{n} \\ 0 \end{bmatrix}, &
      P_{v} & = \begin{bmatrix} 0 & I_{n} \end{bmatrix} P_{\infty}
      \begin{bmatrix} 0 \\ I_{n} \end{bmatrix},\\
      Q_{p} & = \begin{bmatrix} I_{n} & 0 \end{bmatrix} Q_{\infty}
      \begin{bmatrix} I_{n} \\ 0 \end{bmatrix}, &
      Q_{v} & = \begin{bmatrix} 0 & I_{n} \end{bmatrix} Q_{\infty}
      \begin{bmatrix} 0 \\ I_{n} \end{bmatrix}.
  \end{aligned}
\end{align}
Therefor, the extension of the existing second-order balanced truncation
methods to the frequency-limited approach can be done by replacing the
infinite first-order Gramians $P_{\infty}$ and $Q_{\infty}$ 
in~\eqref{eqn:sogramproj} by the first-order frequency-limited
Gramians $P_{\Omega}$ and $Q_{\Omega}$ from~\eqref{eqn:foflgram} corresponding
to the first-order realization~\eqref{eqn:socomp}.
The frequency-limited second-order Gramians are then given by
\begin{align}\label{eqn:soflgramproj}
  \begin{aligned}
    P_{\Omega,p} & = \begin{bmatrix} I_{n} & 0 \end{bmatrix} P_{\Omega}
      \begin{bmatrix} I_{n} \\ 0 \end{bmatrix}, &
      P_{\Omega,v} & = \begin{bmatrix} 0 & I_{n} \end{bmatrix} P_{\Omega}
      \begin{bmatrix} 0 \\ I_{n} \end{bmatrix},\\
    Q_{\Omega,p} & = \begin{bmatrix} I_{n} & 0 \end{bmatrix} Q_{\Omega}
      \begin{bmatrix} I_{n} \\ 0 \end{bmatrix}, &
      Q_{\Omega,v} & = \begin{bmatrix} 0 & I_{n} \end{bmatrix} Q_{\Omega}
      \begin{bmatrix} 0 \\ I_{n} \end{bmatrix},
  \end{aligned}
\end{align}
where $P_{\Omega,p}$ and $P_{\Omega,v}$ are the frequency-limited position
and velocity controllability Gramians, and $J^{\trans}Q_{\Omega,p}J$ and 
$M^{\trans}Q_{\Omega,v}M$ are the frequency-limited position and velocity 
observability Gramians.
Note that $P_{\Omega}$ and $Q_{\Omega}$ are given by~\eqref{eqn:fofllyap}
using the first-order realization~\eqref{eqn:socomp}.
As for the infinite Gramians, one observes that the frequency-limited position
and velocity Gramians are symmetric positive semi-definite.

According to~\cite{morReiS08, morHaiGIetal18, morGawJ90}, we can now
define the corresponding frequency-limited characteristic values as follows.

\begin{definition} \label{def:soflhsv}
  (Second-order frequency-limited characteristic singular values.)\\
  Consider the second-order system~\eqref{eqn:sosys} with
  the first-order realization~\eqref{eqn:socomp} and the frequency range of
  interest $\Omega = -\Omega \subset \mathbb{R}$.
  \begin{enumerate}
    \item The square-roots of the eigenvalues of 
      $P_{\Omega,p}J^{\trans}Q_{\Omega,p}J$ are the \emph{frequency-limited
      position singular values} of~\eqref{eqn:sosys}.
    \item The square-roots of the eigenvalues of 
      $P_{\Omega,p}M^{\trans}Q_{\Omega,v}M$ are the \emph{frequency-limited 
      position-velocity singular values} of~\eqref{eqn:sosys}.
    \item The square-roots of the eigenvalues of 
      $P_{\Omega,v}J^{\trans}Q_{\Omega,p}J$ are the \emph{frequency-limited 
      velocity-po\-si\-tion singular values} of~\eqref{eqn:sosys}.
    \item The square-roots of the eigenvalues of 
      $P_{\Omega,v}M^{\trans}Q_{\Omega,v}M$ are the \emph{frequency-limited 
      velocity singular values} of~\eqref{eqn:sosys}.
  \end{enumerate}
\end{definition}

Following the observations in the first-order frequency-limited case
as well as the second-order balanced truncation method, those
characteristic singular can be interpreted as a measure for the influence
of the corresponding states to the input-output behavior of the system in
the frequency range of interest.
Anyway, there is no energy interpretation as for the first-order balanced 
truncation method.

With~\eqref{eqn:soflgramproj} and the Definition~\ref{def:soflhsv},
the resulting second-order frequency-limited balanced truncation square-root
method is written in Algorithm~\ref{alg:soflbt}.

\begin{algorithm}[t]
  \DontPrintSemicolon
  \caption{Second-Order Frequency-Limited Balanced Truncation Square-Root Method}
  \label{alg:soflbt}
  
  \KwIn{System matrices $M$, $E$, $K$, $B_{u}$, $C_{p}$, $C_{v}$
    from~\eqref{eqn:sosys}, frequency range of interest $\Omega$.}
  \KwOut{Matrices of the reduced-order system $\hM$, $\hE$, $\hK$, $\hB_{u}$
    $\hC_{p}$, $\hC_{v}$.}
    
  Compute Cholesky factorizations of the first-order frequency-limited Gramians 
    by solving~\eqref{eqn:fofllyap}, where the linearization~\eqref{eqn:socomp}
    is used, such that
    $P_{\Omega} = R_{\Omega}R_{\Omega}^{\trans}$,
    $Q_{\Omega} = L_{\Omega}L_{\Omega}^{\trans}$.\;
    
  Follow the steps 2--4 in Algorithm~\ref{alg:sobt}.\;
\end{algorithm}

\begin{remark}
  The second-order frequency-limited balanced truncation method is in general
  not stability preserving.
  Also, the approach from~\cite{morHaiGIetal18} does not necessarily lead to a
  one-sided projection as suggested by the authors and also might not produce
  a stable second-order system in the end.
  Even so, we will discuss approaches that can have stability-preserving
  properties in Section~\ref{sec:modgram}.
\end{remark}


\subsubsection{Second-order time-limited approach}%
\label{sec:sotlbt}

The extension of the time-limited balanced truncation to the second-order
system case was first discussed in~\cite{morHaiGIetal18a}.
As in the previous section, we are generalizing the ideas 
from~\cite{morHaiGIetal18a} to all second-order balanced truncation methods.
In any case, the same idea as for the frequency-limited case is applied here.
That means, we replace the infinite first-order Gramians 
in~\eqref{eqn:sogramproj} by the first-order time-limited Gramians 
from~\eqref{eqn:fotlgram} to get
\begin{align*}
  \begin{aligned}
    P_{T,p} & = \begin{bmatrix} I_{n} & 0 \end{bmatrix} 
      P_{T}\begin{bmatrix} I_{n} \\ 0 \end{bmatrix}, &
      P_{T,v} & = \begin{bmatrix} 0 & I_{n} \end{bmatrix} 
      P_{T} \begin{bmatrix} 0 \\ I_{n} \end{bmatrix},\\
    Q_{T,p} & = \begin{bmatrix} I_{n} & 0 \end{bmatrix} 
      Q_{T}\begin{bmatrix} I_{n} \\ 0 \end{bmatrix}, &
      Q_{T,v} & = \begin{bmatrix} 0 & I_{n} \end{bmatrix} 
      Q_{T}\begin{bmatrix} 0 \\ I_{n} \end{bmatrix},
  \end{aligned}
\end{align*}
where again the first-order realization~\eqref{eqn:socomp} was used.
Following the naming scheme of~\cite{morReiS08}, $P_{T,p}$ and
$P_{T,v}$ are the time-limited position and velocity controllability Gramians,
and $J^{\trans}Q_{T,p}J$ and $M^{\trans}Q_{T,v}M$ the time-limited position and
velocity observability Gramians.
Note that $P_{T}$ and $Q_{T}$ are given 
by~\eqref{eqn:fotllyap} with the first-order realization~\eqref{eqn:socomp}.
As for the infinite Gramians, one observes that the time-limited position
and velocity Gramians are symmetric positive semi-definite.
According to the frequency-limited characteristic singular values, we are 
giving the following definition for the time-limited version.

\begin{definition} \label{def:sotlhsv}
  (Second-order time-limited characteristic singular values.)\\
  Consider the second-order system~\eqref{eqn:sosys} with
  the first-order realization~\eqref{eqn:socomp} and the time range of interest
  $T = [t_{0}, t_{f}]$, $0 \leq t_{0} < t_{f}$.
  \begin{enumerate}
    \item The square-roots of the eigenvalues of $P_{T,p}J^{\trans}Q_{T,p}J$ are
      the \emph{time-limited position singular values} of~\eqref{eqn:sosys}.
    \item The square-roots of the eigenvalues of $P_{T,p}M^{\trans}Q_{T,v}M$ are
      the \emph{time-limited position-velocity singular values}
      of~\eqref{eqn:sosys}.
    \item The square-roots of the eigenvalues of $P_{T,v}J^{\trans}Q_{T,p}J$ are
      the \emph{time-limited velocity-position singular values}
      of~\eqref{eqn:sosys}.
    \item The square-roots of the eigenvalues of $P_{T,v}M^{\trans}Q_{T,v}M$ are
      the \emph{time-limited velocity singular values} of~\eqref{eqn:sosys}.
  \end{enumerate}
\end{definition}

As before, the resulting second-order time-limited balanced truncation methods
can be obtained by replacing the Gramians in Algorithm~\ref{alg:sobt},
which is summarized in Algorithm~\ref{alg:sotlbt}.

\begin{algorithm}[t]
  \DontPrintSemicolon
  \caption{Second-Order Time-Limited Balanced Truncation Square-Root Method}
  \label{alg:sotlbt}
  
  \KwIn{System matrices $M$, $E$, $K$, $B_{u}$, $C_{p}$, $C_{v}$
    from~\eqref{eqn:sosys}, time range of interest $T$.}
  \KwOut{Matrices of the reduced-order system $\hM$, $\hE$, $\hK$, $\hB_{u}$
    $\hC_{p}$, $\hC_{v}$.}
    
  Compute Cholesky factorizations of the first-order time-limited Gramians 
    by solving~\eqref{eqn:fotllyap}, where the linearization~\eqref{eqn:socomp}
    is used, such that
    $P_{T} = R_{T}R_{T}^{\trans}$,
    $Q_{T} = L_{T}L_{T}^{\trans}$.\;
    
  Follow the steps 2--4 in Algorithm~\ref{alg:sobt}.\;
\end{algorithm}

\begin{remark}
  As in the first-order case~\cite{morKue18}, there is no guarantee of
  stability preservation for the second-order time-limited balanced truncation
  methods.
  The method suggested in~\cite{morHaiGIetal18a} only works on the first-order
  case and does not guarantee the preservation of stability for second-order 
  systems in general.
  Approaches that can be more beneficial in terms of preserving stability
  are discussed in Section~\ref{sec:modgram}.
\end{remark}


\section{Numerical methods}%
\label{sec:nummeth}

In this section, we will discuss points concerning the numerical implementation
of the proposed second-order frequency- and time-limited balanced truncation
methods.


\subsection{Matrix equation solvers for large-scale systems}%
\label{sec:mesolvers}

A substantial part of the numerical effort in the computations of the
second-order frequency- and time-limited balanced truncations goes into the
solution of the arising matrix equations~\eqref{eqn:fofllyap} 
and~\eqref{eqn:fotllyap}.
In general it has been shown for the first-order case, that the singular
values of the frequency- and time-limited Gramians are decaying possibly
faster than of the infinite Gramians; see, e.g.,~\cite{morBenKS16}
for the frequency-limited case.
That leads to the natural approximation of the Gramians by low-rank factors, 
e.g.,
\begin{align*}
  \begin{aligned}
    P_{\Omega} & \approx Z_{\Omega}Z_{\Omega}^{\trans}, &
      P_{T} & \approx Z_{T}Z_{T}^{\trans},
  \end{aligned}
\end{align*}
where $Z_{\Omega} \in \mathbb{R}^{N \times \ell_{1}}$, 
$Z_{T} \in \mathbb{R}^{N \times \ell_{2}}$ and
$\ell_{1}, \ell_{2} \ll N$.
Those low-rank factors then replace the Cholesky factors in the balanced
truncation algorithms~\ref{alg:bt}--\ref{alg:sotlbt}.

In the following three sections, we will shortly review existing approaches
for these problems and give comments on existing implementations.


\subsubsection{Quadrature-based approaches}%
\label{sec:quadrature}

A natural approach based on the frequency and time domain integral
representations of the limited Gramians~\eqref{eqn:foflgram}
and~\eqref{eqn:fotlgram} is the use of numerical integration formulas.
As used for example in~\cite{morImrG15, morHaiGIetal18}, the low-rank
factors of the Gramians can be computed by rewriting the full Gramians by
quadrature formulas, e.g.,
\begin{align*}
  P_{\Omega} & = \frac{1}{2\pi}\int\limits_{\Omega}
    (j\omega \cE - \cA)^{-1}\cB\cB^{\trans}(-j\omega \cE - \cA)^{-\trans}
    \mathrm{d}\omega\\
  & \approx \frac{1}{2\pi}\sum\limits_{k = 1}^{\ell} \gamma_{k}
    \{(j\omega_{k}\cE - \cA)^{-1}\cB\cB(-j\omega_{k}\cE - \cA)^{-\trans} +
    (-j\omega_{k}\cE - \cA)^{-1}\cB\cB(j\omega_{k}\cE - \cA)^{-\trans}\},
\end{align*}
where $\gamma_{k}$ are the weights and $\omega_{k}$ the evaluation points
of a fitting quadrature rule, which can be again rewritten for the low-rank
factors by
\begin{align*}
  Z_{\Omega} & = \begin{bmatrix} \real(B_{1}), & \imag(B_{1}), &\ldots
    & \real(B_{\ell}), & \imag(B_{\ell}) \end{bmatrix},
\end{align*}
where $B_k = (j\omega_{k}\cE - \cA)^{-1}\cB$.
Note that this approach becomes unhandy considering the time-limited case,
since there, for each step of the quadrature rule, an approximation of the
matrix exponential has to be computed.

A different approach was suggested in~\cite{morBenKS16}, which writes the
right-hand side of the frequency-limited Lyapunov
equations~\eqref{eqn:fofllyap} as integral expressions, such that
the right-hand side is first approximated and afterwards the large-scale matrix
equation is solved, using one of the approaches in Section~\ref{sec:adi}
or~\ref{sec:krylov}.
In general it is possible to approximate the right-hand sides 
of~\eqref{eqn:fofllyap} and~\eqref{eqn:fotllyap} with matrix functions by
using the general quadrature approach from~\cite{Hig08}.
We are not aware of a stable, available implementation of quadrature-based
matrix equation solvers for the frequency- and time-limited Lyapunov equations
and, therefor, use the following approaches rather than the quadrature-based
methods.


\subsubsection{Low-rank ADI method}%
\label{sec:adi}

The low-rank alternating direction implicit (LR-ADI)~\cite{LiW02, BenLP08}
method is a well established procedure for the solution of large-scale sparse
Lyapunov equations.
Originally developed for the Lyapunov equations corresponding to the
infinite Gramians~\eqref{eqn:folyap}, the LR-ADI produces
low-rank approximations of the form $Z_{\infty, j} = [Z_{\infty, j-1},
\hat{\alpha}_{j}V_{j}]$ by
\begin{align*}
  \begin{aligned}
    V_{j} & = (\cA + \alpha_{j}\cE)^{-1}W_{j-1}, &
      W_{j} & = W_{j-1} - 2\real(\alpha_{j})V_{j},
  \end{aligned}
\end{align*}
where $\hat{\alpha}_{j} = \sqrt{-2\real{\alpha_{j}}}$, $W_{0} = \cB$;
see~\cite{morBenKS16, BenKS13a, BenKS14b} for more details on this method.

The right-hand sides of the limited Lyapunov 
equations~\eqref{eqn:fofllyap},~\eqref{eqn:fotllyap} can be rewritten as
\begin{align} \label{eqn:indefrhs}
  \begin{aligned}
    \cB_{\Omega}\cB^{\trans} + \cB\cB_{\Omega}^{\trans} & =
      \widetilde{\cB}
      \begin{bmatrix} 0 & I_{m} \\ I_{m} & 0 \end{bmatrix} 
      \widetilde{\cB}^{\trans}, & 
      \cC_{\Omega}^{\trans}\cC + \cC^{\trans}\cC_{\Omega} & =
      \widetilde{\cC}^{\trans}
      \begin{bmatrix} 0 & I_{p} \\ I_{p} & 0 \end{bmatrix}
      \widetilde{\cC},\\
    \cB_{t_{0}}\cB_{t_{0}}^{\trans} - \cB_{t_{f}}\cB_{t_{f}}^{\trans} & =
      \breve{\cB}\begin{bmatrix} I_{m} & 0 \\ 0 & -I_{m} \end{bmatrix}
      \breve{\cB}^{\trans}, &
      \cC_{t_{0}}^{\trans}\cC_{t_{0}} - \cC_{t_{f}}^{\trans}\cC_{t_{f}} & =
      \breve{\cC}^{\trans}
      \begin{bmatrix} I_{p} & 0 \\ 0 & -I_{p} \end{bmatrix}\breve{\cC}
  \end{aligned}
\end{align}
with $\widetilde{\cB} = [\cB_{\Omega},\cB]$,
$\widetilde{\cC}^{\trans} = [\cC_{\Omega}^{\trans},\cC^{\trans}]$,
$\breve{\cB} = [\cB_{t_{0}},\cB_{t_{f}}]$ and
$\breve{\cC}^{\trans} = [\cC_{t_{0}}^{\trans}, \cC_{t_{f}}^{\trans}]$,
which shows that the right-hand side matrices are indefinite.
The LR-ADI method can be extended to this case by using an 
$LDL^{\trans}$-factorization for the right-hand side as well as for the 
solution~\cite{LanMS15}.
Note that for applying this method for the solution of the large-scale matrix 
equations, an approximation of the matrix functions in the right-hand sides is 
needed beforehand.
It was noted in~\cite{morBenKS16}, that the information used for the 
approximation of the matrix functions cannot be used in the LR-ADI method.
A stable version of the LR-ADI method in the low-rank and $LDL^{\trans}$
formats is implemented in~\cite{SaaKB19-mmess-2.0}.
We will use this implementation in case the methods, described in the following
section, are failing to converge for the solution of the matrix equation but
give approximations to the function right hand-sides.

  
\subsubsection{Projection methods}%
\label{sec:krylov}

An approach that can be used to approximate the matrix functions in the
right-hand sides of the limited Lyapunov equations, as well as to solve the
large-scale matrix equations at the same time, is given by projection-based
methods.
Here, low-dimensional subspaces $\cV_{k} = \range(V_{k})$ are used to obtain
the low-rank solutions as, e.g., $P_{\Omega} \approx V_{k}
\check{P}_{\Omega} V_{k}^{\trans}$, where $\check{P}_{\Omega}$ is the solution
of the projected Lyapunov equation
\begin{align} \label{eqn:projlyap}
  \cT_{k} \check{P}_{\Omega} + \check{P}_{\Omega} \cT_{k}^{\trans} +
    \check{\cB}_{\Omega}\check{\cB}^{\trans} + 
    \check{\cB}\check{\cB}_{\Omega}^{\trans} & = 0,
\end{align}
$\cT_{k} = V_{k}^{\trans} \cE^{-1}\cA V_{k}$,
$\check{\cB}_{\Omega} = V_{k}^{\trans} \cE^{-1} \cB_{\Omega}$ and
$\check{\cB} = V_{k}^{\trans} \cE^{-1} \cB$ are the projected matrices 
of the frequency-limited controllability Lyapunov
equation~\eqref{eqn:fofllyap}.
The equation~\eqref{eqn:projlyap} is now small and dense and can be solved
using established dense solvers.
As one can observe, this method gives also the opportunity to approximate the 
matrix function right-hand side by the low-dimensional subspace
$\cV_{k}$, for which one can also use dense computation methods~\cite{Hig08}.

Usually, the low-dimensional subspace $\cV_{k}$ is constructed
as standard~\cite{morJaiK94}, extended~\cite{Sim07} or rational Krylov
subspace~\cite{DruS11}, all of which can be easily computed for
large-scale sparse systems.
The implementation of the limited balanced truncation
methods for second-order systems~\cite{morBenW20}, we provide, is also
based on rational Krylov subspaces.
We refer the reader to~\cite[Algorithm 4.1]{morBenKS16} for the underlying idea
of the implementation.

A drawback of the projection-based approach, especially for second-order
systems, is that the projected system matrices $\cT_{k}$ are not necessarily
c-stable, even if the original first-order realization of the second-order
system was.
In fact, the quality and performance of the projection-based solvers strongly
depend on the chosen first-order realization.
Therefor, we are going to use the so-called strictly dissipative realization of
second-order systems~\cite{PanWL12} in our computations.
Assuming $M, E, K$ to be symmetric positive definite, the second-order
system~\eqref{eqn:sosys} can be described by a first-order realization using
the following matrices
\begin{align} \label{eqn:stricdissrel}
  \begin{aligned}
    \cE & = \begin{bmatrix} K & \gamma M \\ \gamma M & M \end{bmatrix}, &
    \cA & = \begin{bmatrix} -\gamma K & K - \gamma E \\ -K & -E + \gamma M
      \end{bmatrix}, &
    \cB & = \begin{bmatrix} \gamma B_{u} \\ B_{u} \end{bmatrix}, &
    \cC & = \begin{bmatrix} C_{p} & C_{v} \end{bmatrix},
  \end{aligned}
\end{align}
with the parameter $0 < \gamma < \lambda_{\min}(E(M +
\frac{1}{4}E K^{-1} E)^{-1})$.
The advantage of this realization is that $\cE$ is symmetric positive definite
and $\cA + \cA^{\trans}$ symmetric negative definite.
Following that, projection methods can preserve the stability in the
projected matrices $\cT_{k}$ if the computations are made on the corresponding
standard state-space realization, obtained by a symmetric state-space
transformation using the Cholesky factors $\cE = \cL\cL^{\trans}$, i.e.,
the algorithms work implicitly on a realization of the form
\begin{align*}
  \begin{aligned}
    \dot{\tq}(t) & = \cL^{-1}\cA\cL^{-\trans} \tq(t) + \cL^{-1}\cB u(t),\\
    y(t) & = \cC\cL^{-\trans} \tq(t).
  \end{aligned}
\end{align*}

\begin{remark}
  Note that the realization~\eqref{eqn:stricdissrel} is computationally more
  involved than the classical first companion form~\eqref{eqn:socomp} or its
  second companion form, since it is not possible to make use of occurring
  zeros in the block structure.
\end{remark}

Also, by changing the first-order realization to~\eqref{eqn:stricdissrel},
the computed Gramians change compared to the definition of the second-order
balancing methods.
Therefor, let $\tP$ and $\tQ$ be Gramians computed for the strictly
dissipative first-order realization~\eqref{eqn:stricdissrel} and $P$ and $Q$
be the Gramians from the first companion form realization~\eqref{eqn:socomp}.
Then it holds
\begin{align*}
  \begin{aligned}
    P & = \tP && \text{and} & Q & = T^{\trans} \tQ T,
  \end{aligned}
\end{align*}
with the transformation matrix
\begin{align*}
  T = & \begin{bmatrix} K & \gamma I_{n} \\ \gamma M & I_{n} \end{bmatrix}.
\end{align*}
That means we can use the strictly dissipative
realization~\eqref{eqn:stricdissrel} for the solution of the matrix equations
and for the balancing procedure just perform the easy back transformation of the
observability factor.


\subsection{Stabilization and acceleration by %
  \texorpdfstring{$\alpha$}{alpha}-shifts}%
\label{sec:alphashifts}

So far, it was always assumed that the second-order system~\eqref{eqn:sosys} is
c-stable.
But in practice, the eigenvalues of $\lambda^{2}M + \lambda E + K$ can be
very close to the imaginary axis or even on the axis, e.g., in the case
of marginal stability.
This makes the usage of the model reduction methods and matrix equation
solvers very difficult.
A strategy to overcome those problems has been proposed in,
e.g.,~\cite{morFreRM08}.
There, a shift in the frequency domain was used to move the spectrum of the
pencil $\lambda\cE - \cA$, which had eigenvalues at zero, away from the
imaginary axis to compute the system Gramians.
This approach cannot be used the same way for the first-order 
realizations~\eqref{eqn:socomp} or~\eqref{eqn:stricdissrel} of second-order
systems since it destroys the block structure one can exploit in the
numerical implementations of the solvers or rather the block structure that
is used for the second-order balancing approaches.
Therefore, we will transfer the concept of $\alpha$-shifts to the case of
second-order systems.

Let $\alpha \in \mathbb{R}_{> 0}$ be a real, strictly positive shift and
consider the second-order differential equations in the frequency-domain
\begin{subequations}
\begin{align} \label{eqn:sosysfreq1}
  (s^{2}M + sE + K)X(s) & = B_{u}U(s),\\ \label{eqn:sosysfreq2}
  Y(s) & = (sC_{v} + C_{p})X(s),
\end{align}
\end{subequations}
where $U(s), X(s), Y(s)$ are the Laplace transforms of the corresponding time
domain functions and $s \in \mathbb{C}$ the Laplace variable.
Now let $s = \rho + \alpha$, with $\rho \in \mathbb{C}$ a shifted Laplace
variable.
Then the equation~\eqref{eqn:sosysfreq1} turns into
\begin{align*}
  ((\rho + \alpha)^{2}M + (\rho + \alpha)E + K)X(s)
    & = (\rho^{2}M + 2\alpha\rho M + \alpha^{2}M
    + \rho D + \alpha E + K)X(s)\\
  & = (\rho^{2}M + \rho(E + 2\alpha M) + (K + \alpha E + \alpha^{2}M))X(s)\\
  & = (\rho^{2}M + \rho\tE + \tK)X(s)\\
  & = B_{u}U(s),
\end{align*}
with $\tE = E + 2\alpha M$ and $\tK = K + \alpha E + \alpha^2 M$.
Also, the second equation~\eqref{eqn:sosysfreq2} can be rewritten as
\begin{align*}
  Y(s) & = ((\rho + \alpha)C_{v} + C_{p})X(s)\\
  & = (\rho C_{v} + (C_{p} + \alpha C_{v}))X(s)\\
  & = (\rho C_{v} + \tC_{p})X(s),
\end{align*}
where $\tC_{p} = C_{p} + \alpha C_{v}$.
Now, the new system described by $(M, \tE, \tK, B_{u}, \tC_{p}, C_{v})$ is used
for the computation of the reduced-order projection matrices
$W, T \in \mathbb{R}^{n \times r}$.
Then, the projected system $(\hM, \widehat{\tE}, 
\widehat{\tK}, \hB_{u}, \widehat{\tC}_{p}, C_{v})$
yields the following relations
\begin{align*}
  \begin{aligned}
    \widehat{\tE} & = \hE + 2\alpha \hM, &
    \widehat{\tK} & = \hK + \alpha \hE + \alpha^{2} \hM, &
    \widehat{\tC}_{p} & = \hC_{p} + \alpha \hC_{v},
  \end{aligned}
\end{align*}
where $\hE = W^{\trans} E T$, $\hK = W^{\trans} K T$ and $\hC_{p} = C_{p} T$
are the transformed non-shifted matrices.
Now, we consider the transformed system again in the frequency domain with the
Laplace variable $\rho$ and using the back-substitution $\rho = s - \alpha$,
such that
\begin{align*}
  \begin{aligned}
    \rho^{2}\hM + \rho \widehat{\tE} + \widehat{\tK}
      & = s^{2}\hM + s\hE + \hK
      & \text{and} &&
      \rho \hC_{v} + \widehat{\tC}_{p} & = s\hC_{v} + \hC_{p}.
  \end{aligned}
\end{align*}
The back-substitution gives the resulting reduced-order model
$(\hM, \hE, \hK, \hB_{u}, \hC_{p}, \hC_{v})$.
The $\alpha$-shift strategy can be interpreted as a structured perturbation
in the frequency domain during the computations.
Experiments have shown that such an approach works fine for $\alpha$ small
enough.
It has to be noted that there are no theoretical results on the influence
of the chosen $\alpha$ concerning the quality of the reduced-order model or
properties like stability preservation and error bounds.

\begin{remark}
  The $\alpha$-shift approach can also be used either to improve the
  conditioning of the used matrix equation solvers by improving the condition
  number of the shifted linear systems solving with $(\sigma^{2}M + \sigma \tE
  + \tK)$, or to improve the convergence of those solvers by pushing the
  eigenvalues of $\lambda^{2}M + \lambda \tE + \tK$ further away from the
  imaginary axis.
\end{remark}


\subsection{Two-step hybrid methods}%
\label{sec:hybridmor}

The idea of two-step (or hybrid) model reduction methods has been used for
quite some time in different
applications~\cite{morLehE07, morFehE10, morWolPL13}.
In general, two-step methods are based on the division of the model reduction 
process into two phases.
First a pre-reduction, which can be easily computed and gives a very
accurate approximation for the system's behavior.
The model resulting from the pre-reduction is usually of medium-scale
dimensions, on which the second reduction step by a more sophisticated model
reduction method is applied.
This procedure has the advantage that there is no necessity of applying
difficult approximation methods for the large-scale matrix equations
arising in the balancing related approaches.
Instead, the exact methods can be used on the, usually dense,
pre-reduced system.

In order to have a structure-preserving pre-reduction method, we suggest
the use of interpolation by rational Krylov 
subspaces~\cite{morBeaG05, morSal05, morSalL06}.
This has been shown to be equivalent to the use of shift-based approximation
methods for the large-scale matrix equations in Section~\ref{sec:mesolvers}; 
see~\cite{morWolPL13}.
The second-order rational Krylov subspaces are generated as
\begin{align*}
  \cV & = \range\left((s_{1}^{2}M + s_{1}E + K)^{-1}B_{u},\ldots, 
    \prod\limits_{k = 1}^{\ell}(s_{k}^{2}M + s_{k}E + K)^{-1}B_{u}\right),\\
  \cU & = \range\left((s_{1}^{2}M + s_{1}E + K)^{-\herm}
    (C_{p} + s_{1}C_{v})^{\herm},\ldots, 
      \prod\limits_{k = 1}^{\ell}(s_{k}^{2}M + s_{k}E + K)^{-\herm}
      (C_{p} + s_{k}C_{v})^{\herm}\right),
\end{align*}
with $s_{k} \in \mathbb{C}$, $k = 1, \ldots, \ell$, chosen interpolation points.
Let $V$ and $U$ be Hermitian bases of the same size such that
$\cV \subset \range(V)$ and $\cU \subset \range(U)$, respectively, the
pre-reduced model is then generated by
\begin{align*}
  \begin{aligned}
    M_{\pre} & = U^{\herm}MV, & E_{\pre} & = U^{\herm}EV, & 
      K_{\pre} & = U^{\herm}KV,\\
    B_{u, \pre} & = U^{\herm}B_{u}, & C_{p, \pre} & = C_{p} V, &
      C_{v, \pre} & = C_{v} V.
  \end{aligned}
\end{align*}
For preservation of stability and the realness of the system matrices, we
choose the interpolation points to appear in complex conjugate pairs $s_{k}$ and
$\overline{s}_{k}$, and replace one of the projection matrices by $U = V$.

The choice of points $s_{k}$ is crucial for the quality of the pre-reduced
model.
While there are strategies for an adaptive or optimal choice of $s_{k}$,
we suggest a simple oversampling on the imaginary axis, which
is usually enough as a global pre-reduced model.

\begin{remark}
  For the frequency-limited case, a natural choice for the interpolation 
  points would be to take $j\Omega$ instead of aiming for a global
  approximation.
  In this case, the resulting frequency-limited balanced truncation
  will very likely not give the same results as the large-scale approach.
  This observation comes from the fact, that the frequency-limited balanced
  truncation still takes information about the complete system structure into 
  account and the pre-reduced system can be completely different from the
  original one, if only a local pre-reduction is performed.
\end{remark}

Due to the required accuracy of the pre-reduced model, the dimension of it can
be still very large.
Therefore, we suggest an efficient iterative solver for the Lyapunov equations
appearing in the second reduction step.
In general, we consider the following stable Lyapunov equations
\begin{align} \label{eqn:indeflyap}
  \begin{aligned}
    \cA X_{1} \cE^{\trans} + \cE X_{1} \cA^{\trans} + \cB \cQ \cB^{\trans}
      & = 0,\\
    \cA^{\trans} X_{2} \cE + \cE^{\trans} X_{2} \cA + \cC^{\trans} \cR \cC
      & = 0,
  \end{aligned}
\end{align}
where $\cQ \in \mathbb{R}^{m \times m}$ and $\cR \in \mathbb{R}^{p \times p}$
are symmetric and possibly indefinite.
The solution of~\eqref{eqn:indeflyap} can then be factored in the same way
as the right-hand sides, i.e., $X_{1} = Z_{1}Y_{1}Z_{1}^{\trans}$ and
$X_{2} = Z_{2}Y_{2}Z_{2}^{\trans}$, where $Y_{1}$ and $Y_{2}$ are also
symmetric matrices.
For efficiently computing the solutions of~\eqref{eqn:indeflyap},
we extend the dual sign function iteration method from~\cite{BenCQ98a}
for the $LDL^{\trans}$-factorization of the solutions.
As a result, we get a sign function iteration, that solves both Lyapunov 
equations with symmetric indefinite right hand-sides~\eqref{eqn:indeflyap} at
the same time; see  Algorithm~\ref{alg:lyapsgnldl}.

\begin{algorithm}[tb]
  \DontPrintSemicolon
  \caption{$LDL^{\trans}$-Factored Sign Function Dual Lyapunov Equation Solver}
  \label{alg:lyapsgnldl}
  
  \KwIn{$\cA$, $\cB$, $\cC$, $\cE$, $\cQ$, $\cR$ from~\eqref{eqn:indeflyap},
    tolerance $\tau$.}
  \KwOut{$Z_{1}$, $Y_{1}$, $Z_{2}$, $Y_{2}$ -- solution factors of
   ~\eqref{eqn:indeflyap}.}
  
  Set $A_{1} = \cA$, $B_{1} = \cB$, $Q_{1} = \cQ$, $C_{1} = \cC$,
    $R_{1} = \cR$, $k = 1$.\;
  \While{$\lVert A_{k} + \cE \rVert > \tau \lVert E \rVert$}{
    Compute the scaling factor for convergence acceleration
      \begin{align*}
        c_{k} & = \sqrt{\frac{\lVert A_{k} \rVert_{F}}
          {\lVert \cE A_{k}^{-1} \cE \rVert_{F}}}.
      \end{align*}\vspace{-\baselineskip}\;
    Compute the next iterates of the solution factors
    \begin{align*}
      \begin{aligned}
        B_{k+1} & = \begin{bmatrix} B_{k}, & \cE A_{k}^{-1} B_{k}\end{bmatrix},&
          Q_{k+1} & = \begin{bmatrix} \frac{1}{2c_{k}} Q_{k} & \\ & 
          \frac{c_{k}}{2} Q_{k} \end{bmatrix},\\
        C_{k+1} & = \begin{bmatrix}C_{k} \\ A_{k}^{-1}\cE C_{k}\end{bmatrix},&
          R_{k+1} & = \begin{bmatrix} \frac{1}{2c_{k}} R_{k} & \\ & 
          \frac{c_{k}}{2} R_{k} \end{bmatrix}.
      \end{aligned}
    \end{align*}\vspace{-\baselineskip}\;
    Compute the next iteration matrix
    \begin{align*}
      A_{k+1} = \frac{1}{2c_{k}}A_{k} + \frac{c_{k}}{2}\cE A_{k}^{-1}\cE
    \end{align*}\vspace{-\baselineskip}\;
    Set $k = k + 1$.\;
  }
  Construct the solution factors
  \begin{align*}
    \begin{aligned}
      Z_{1} & = \frac{1}{\sqrt{2}} \cE^{-1} B_{k}, & 
        Y_{1} & = Q_{k}, &
        Z_{2} & = \frac{1}{\sqrt{2}} \cE^{-\trans}C_{k}^{\trans}, &
        Y_{2} & = R_{k}.
    \end{aligned}
  \end{align*}\vspace{-\baselineskip}\;
\end{algorithm}

The implementation of Algorithm~\ref{alg:lyapsgnldl}  as well as dense versions
of the second-order frequency- and time-limited balanced truncation methods can 
be found in~\cite{morBenW19b}.

\begin{remark}
  In Step 4 of Algorithm~\ref{alg:lyapsgnldl}, the memory requirements and
  operations are doubling in every step due to the extension of the solution
  factors.
  It is suggested to do $LDL^{\trans}$ column and row compressions at that
  point to keep the size of the factors small.
\end{remark}


\subsection{Modified Gramian approach}%
\label{sec:modgram}

A drawback of the frequency- and time-limited balanced truncation methods
is the loss of stability preservation.
For the first-order system case, there are different modifications of the
methods to regain the preservation of stability, e.g.,
the replacement of one of the limited Gramians by the infinite 
Gramian~\cite{morImrG15, morHaiGIetal17}.

A different technique, proposed in~\cite{morGugA04}, is the modified Gramian
approach.
Therefor, the indefinite right-hand sides~\eqref{eqn:indefrhs} are replaced
by definite ones.
Using eigenvalue decompositions, the right-hand sides can be rewritten as
\begin{align*}
  \begin{aligned}
    \cB_{\Omega}\cB^{\trans} + \cB\cB_{\Omega}^{\trans} & = 
      U_{\cB,\Omega}S_{\cB,\Omega}U_{\cB,\Omega}^{\trans}, &
      \cC_{\Omega}^{\trans}\cC + \cC^{\trans}\cC_{\Omega} & =
      U_{\cC,\Omega}S_{\cC,\Omega}U_{\cC,\Omega}^{\trans}\\
    \cB_{t_{0}}\cB_{t_{0}}^{\trans} - \cB_{t_{f}}\cB_{t_{f}}^{\trans} & = 
      U_{\cB,T}S_{\cB,T}
      U_{\cB,T}^{\trans}, &
      \cC_{t_{0}}^{\trans}\cC_{t_{0}} - \cC_{t_{f}}^{\trans}\cC_{t_{f}} & =
      U_{\cC,T}S_{\cC,T}U_{\cC,T}^{\trans},
  \end{aligned}
\end{align*}
where $U_{\cB,\Omega}$, $U_{\cC,\Omega}$, $U_{\cB,T}$,
$U_{\cC,T}$ are orthogonal and
\begin{align*}
  \begin{aligned}
    S_{\cB,\Omega} & = \diag(\eta^{\cB}_{1}, \ldots, \eta^{\cB}_{2m},
      0, \ldots, 0), &
      S_{\cC,\Omega} & = \diag(\eta^{\cC}_{1}, \ldots, \eta^{\cC}_{2p},
      0, \ldots, 0),\\
    S_{\cB,T} & = \diag(\mu^{\cB}_{1}, \ldots, \mu^{\cB}_{2m},
      0, \ldots, 0), &
      S_{\cC,T} & = \diag(\mu^{\cC}_{1}, \ldots, \mu^{\cC}_{2p},
      0, \ldots, 0).
  \end{aligned}
\end{align*}
Let $U_{\cB,\Omega,1}$, $U_{\cC,\Omega,1}$, $U_{\cB,T,1}$,
$U_{\cC,T,1}$
be the parts of the orthogonal matrices, corresponding to the possible
non-zero eigenvalues.
The modified frequency- and time-limited Gramians are then given as the 
solutions of the following Lyapunov equations
\begin{align*}
  \begin{aligned}
    \cA P_{\Omega}^{\mo} \cE^{\trans} + \cE P_{\Omega}^{\mo} \cA^{\trans}
      + \cB_{\Omega}^{\mo}\left(\cB_{\Omega}^{\mo}\right)^{\trans} & = 0,\\
    \cA^{\trans} Q_{\Omega}^{\mo} \cE + \cE^{\trans} Q_{\Omega}^{\mo} \cA
      + \left(\cC_{\Omega}^{\mo}\right)^{\trans}\cC_{\Omega}^{\mo} & = 0,\\
    \cA P_{T}^{\mo} \cE^{\trans} + \cE P_{T}^{\mo} 
      \cA^{\trans} + \cB_{T}^{\mo}
      \left(\cB_{T}^{\mo}\right)^{\trans} & = 0,\\
    \cA^{\trans} Q_{T}^{\mo} \cE + \cE^{\trans}
      Q_{T}^{\mo} \cA
      + \left(\cC_{T}^{\mo}\right)^{\trans}
      \cC_{T}^{\mo} & = 0,
  \end{aligned}
\end{align*}
with
\begin{align*}
  \begin{aligned}
    \cB_{\Omega}^{\mo} & = U_{\cB,\Omega,1}\diag(\lvert\eta^{\cB}_{1}\rvert, 
      \ldots, \lvert\eta^{\cB}_{2m}\rvert), &
      \cC_{\Omega}^{\mo} & = \diag(\lvert\eta^{\cC}_{1}\rvert, \ldots, 
      \lvert\eta^{\cC}_{2p}\rvert) U_{\cC,\Omega,1}^{\trans},\\
    \cB_{T}^{\mo} & = U_{\cB,T, 1} 
      \diag(\lvert\mu^{\cB}_{1}\rvert, \ldots, \lvert\mu^{\cB}_{2m}\rvert), &
      \cC_{T}^{\mo} & = \diag(\lvert\mu^{\cC}_{1}\rvert,
      \ldots, \lvert\mu^{\cC}_{2p}\rvert)U_{\cC,T,1}^{\trans}.
  \end{aligned}
\end{align*}
Using those modified Gramians for the limited balanced truncation methods also
preserves the stability in the reduced-order models in the first-order case.
There also exists an $\mathcal{H}_{\inf}$ error bound for the modified 
frequency-limited balanced truncation for first-order
systems~\cite{morBenKS16}.
Note that the limited Gramians can also be easily computed using the
projection-based matrix equation solvers with only minor changes in the
algorithm~\cite{morBenKS16, morKue18}.

\begin{remark}
  Neither the replacement of limited Gramians by the infinite ones nor the
  modified Gramian approaches are guaranteed to preserve the stability in the
  reduced-order model when it comes to the second-order case.
  The stability preserving methods in~\cite{morHaiGIetal18, morHaiGIetal18a}
  are just based on the assumption, that the same procedure as in the 
  first-order case also works for second-order systems.
  This is not the case, since already the classical second-order balanced 
  truncation methods are in general not stability preserving~\cite{morReiS08}.
\end{remark}

\begin{remark}
  Also, it has been mentioned and shown by numerical examples
  in~\cite{morBenKS16, morKue18} that the modified Gramian approach usually
  does not pay off since the quality of the reduced-order models is often the
  same as for the global approaches, i.e., the local approximation property of
  the limited balanced truncation methods gets lost.
\end{remark}


\begin{figure}[tb]
  \begin{center}
    \tikzexternalenable
    \tikzsetnextfilename{damped_mass_spring}
    \input{graphics/damped_mass_spring.tikz}
    \tikzexternaldisable
    \vspace{-1.5\baselineskip}
  \end{center}
  \caption{Setup of the single chain oscillator.}
  \label{fig:sco}
\end{figure}
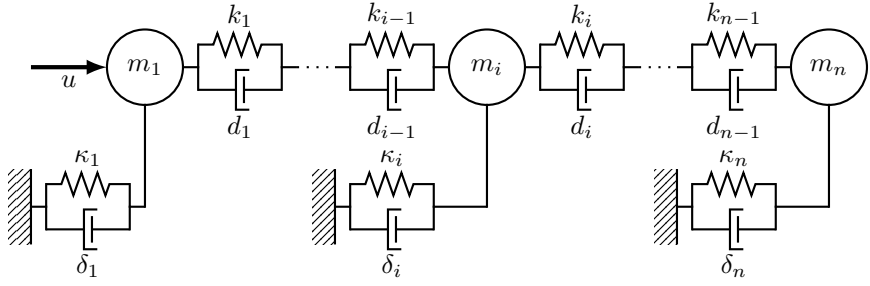

\section{Numerical examples}
\label{sec:examples}

In the following, some mechanical systems of second-order form from the
literature have been chosen as benchmark examples.
The experiments reported here have been executed on machines with 2 Intel(R)
Xeon(R) Silver 4110 CPU processors running at 2.10GHz and equipped with either
192 GB or 384 GB total main memory.
The computers are running on CentOS Linux release 7.5.1804 (Core) and using
MATLAB 9.4.0.813654 (R2018a).
For the computations, the following software has been used:
\begin{itemize}
  \item MORLAB version 5.0~\cite{morBenW19b}, for all evaluations in the
    frequency and time domain, the generation of the pictures and the dense
    implementations of the limited model reduction methods used in the two-step
    approach,
  \item the limited balanced truncation for large-scale sparse second-order
    systems code package~\cite{morBenW20}, for the
    computations of the full-order limited Gramians and the implementation of
    the balancing formulas from Table~\ref{tab:sobt},
  \item the M-M.E.S.S. library version 2.0~\cite{SaaKB19-mmess-2.0}, for
    computing the full Gramians with already approximated right hand-sides.
\end{itemize}
In general, we used the projection-based methods from~\cite{morBenW20} to
approximate the right hand-sides and the Gramians.
But in case that the Gramians did not converge, we used the computed
approximation of the right hand-sides from the projection methods in the ADI
method from~\cite{SaaKB19-mmess-2.0} to compute a solution to the matrix
equation.

For the presentation of the results, the following error measures have been
used.
In the frequency domain, the point-wise absolute errors in the plots are
computed as $\lVert H(j\omega) - \hH(j\omega) \rVert_{2}$ for the frequency
points $\omega \in \mathbb{R}$ and the point-wise relative error as
$\frac{\lVert H(j\omega) - \hH(j\omega) \rVert_{2}}%
{\lVert H(j\omega)\rVert_{2}}$.
The corresponding error tables show as global errors the maximum
value of the point-wise errors in the plotted frequency region, i.e.,
\begin{align*}
  \begin{aligned}
    \max\limits_{\omega \in [\omega_{\min},\omega_{\max}]}
      {\lVert H(j\omega) - \hH(j\omega) \rVert_{2}} && \text{and} &
    \max\limits_{\omega \in [\omega_{\min},\omega_{\max}]}
      {\frac{\lVert H(j\omega) - \hH(j\omega) \rVert_{2}}%
      {\lVert H(j\omega)\rVert_{2}}},
  \end{aligned}
\end{align*}
where $[\omega_{\min},\omega_{\max}]$ is the frequency region as shown in the
plots.
The local errors are then the maximum values in the frequency range of interest.

In the time domain, the errors are also point-wise evaluated.
The plots show $\lVert y(t) - \hy(t) \rVert_{2}$ with $t \in \mathbb{R}$ as
absolute errors and $\frac{\lVert y(t) - \hy(t) \rVert_{2}}%
{\lVert y(t) \rVert_{2}}$ for the relative errors.
The corresponding error tables show again the maximum point-wise error values
\begin{align*}
  \begin{aligned}
    \max\limits_{t \in [t_{\min}, t_{\max}]}
      {\lVert y(t) - \hy(t) \rVert_{2}} && \text{and} &
    \max\limits_{t \in [t_{\min}, t_{\max}]}
      {\frac{\lVert y(t) - \hy(t) \rVert_{2}}{\lVert y(t) \rVert_{2}}},
  \end{aligned}
\end{align*}
where $[t_{\min}, t_{\max}]$ is the time frame as shown in the plots or rather
the local time range $[t_{0}, t_{f}]$ chosen for the time-limited methods.

As criterion for the computed approximation order, the characteristic values
from Definition~\ref{def:soflhsv} and~\ref{def:sotlhsv} have been used.
Therefore, we truncated all states corresponding to the singular values that
in sum were smaller than the largest singular values multiplied with the
tolerance $10^{-4}$, i.e.,
\begin{align*}
  10^{-4}\sigma_{1} \geq \sum\limits_{k = r+1}^{n_{\min}}\sigma_{k}.
\end{align*}


\subsection{Single chain oscillator}

\begin{figure}[tb]
  \begin{center}
    \tikzexternalenable
    \tikzsetnextfilename{singlechain_soflbt_a0}
    \input{graphics/singlechain_soflbt_a0.tikz}
    \tikzexternaldisable
    \vspace{-2\baselineskip}
  \end{center}
  \caption{Frequency-limited ROMs for the single chain oscillator
    (full-order Gramians).}
  \label{fig:sco_fl}
\end{figure}
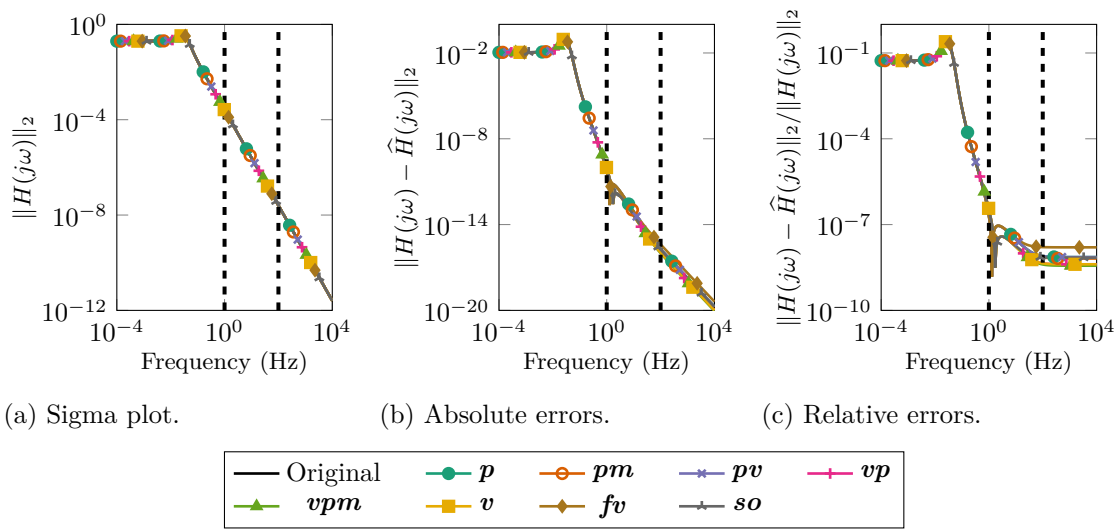

As first example, we consider the single chain oscillator benchmark
from~\cite{morMehS05}, where we removed the holonomic constraint to get a
mechanical system without algebraic parts.
Figure~\ref{fig:sco} shows the basic setup of the system, where the parameters
are chosen as in~\cite{morMehS05}, i.e. in our experiments we have
\begin{align*}
  m_{1} & = \ldots = m_{n} = 100,\\
  k_{1} & = \ldots = k_{n-1} = \kappa_{2} = \ldots = \kappa_{n-1} = 2,\\
  d_{1} & = \ldots = d_{n-1} = \delta_{2} = \ldots = \delta_{n-1} = 5,
\end{align*}
and $\kappa_{1} = \kappa_{n} = 4$, $\delta_{1} = \delta_{n} = 10$.
The input and output matrices are chosen to be $B_{u} = e_{1}$ and
$C_{p} = [e_{1}, e_{2}, e_{n-1}]^{\trans}$, where $e_{i}$ denotes the $i$-th
column of the identity matrix $I_{n}$.
Also, we have chosen $n = 12\,000$ masses for the system.
This system doe not have any velocity outputs $C_{v}$.

\subsubsection{Frequency domain}

The frequency range of interest in this example is chosen, just for
demonstration reasons, to be between $1$ and $100$\,Hz.
The computations have been done with no $\alpha$-shift ($\alpha = 0$).
In Figure~\ref{fig:sco_fl}, the resulting reduced-order models (ROMs) can be
seen in terms of their transfer functions (a), the point-wise absolute error (b)
and point-wise relative error (c).
The frequency range of interest is marked as the area between the dashed
vertical lines.
Table~\ref{tab:sco_fl} gives an overview for all applied second-order
frequency-limited.
It can be noted that all computed ROMs are of order $2$,
stable and have absolute and relative errors in the same order of magnitude.
Also we note that as wanted, the errors in the frequency range of interest are
significantly smaller than in the overall considered frequency region.
For the two-step approach, we used, on the one hand, a logarithmically
equidistant sampling of $200$ frequency points in the frequency region of
interest and, on the other, for a global approximation logarithmically
equidistant points between $10^{-4}$ and $10^{4}$\,Hz.
After a rank truncation of the orthogonalized basis, the intermediate
ROMs had the dimension $100$.
Since no significant differences between the full-order Gramian and two-step
approaches could be seen, we refer the reader also to Figure~\ref{fig:sco_fl}
and Table~\ref{tab:sco_fl} for the results.

\subsubsection{Time domain}

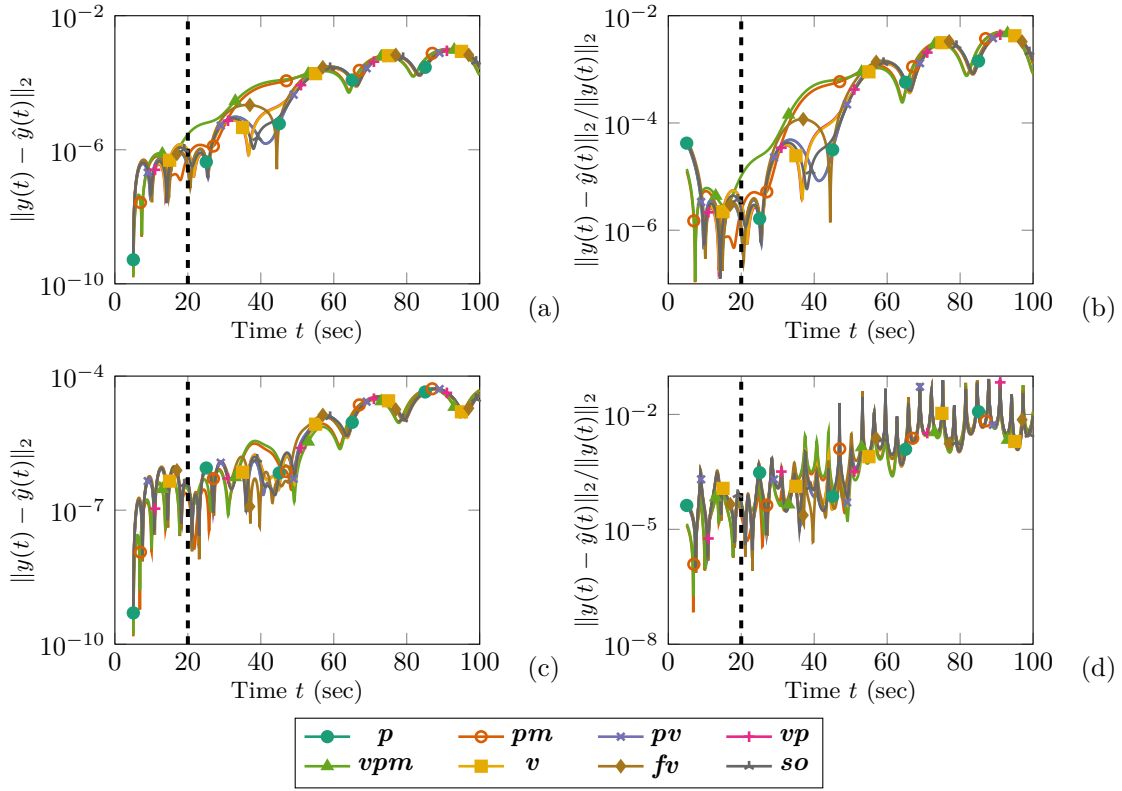
\begin{figure}[tb]
  \begin{center}
    \tikzexternalenable
    \tikzsetnextfilename{singlechain_sotlbt_a0}
    \input{graphics/singlechain_sotlbt_a0.tikz}
    \tikzexternaldisable
    \vspace{-2\baselineskip}
  \end{center}
  \caption{Absolute and relative errors of time-limited ROMs for the single
    chain oscillator with inputs $u_{\step}$ (a), (b) and
    $u_{\sin}$ (c), (d) (full-order Gramians).}
  \label{fig:sco_tl}
\end{figure}

In the time domain, we apply two different input signals to test our ROMs
\begin{align} \label{eqn:input_sco}
  \begin{aligned}
    u_{\step}(t) & = \delta(t - 5) && \text{and} &
      u_{\sin}(t) & = \sin(t)\delta(t - 5),
  \end{aligned}
\end{align}
for $t \in [0, 100]$ and $\delta(t)$ the Heaviside function.
As time range of interest, $[0, 20]$ has been chosen.

While Figure~\ref{fig:sco_tl} shows the results for the time-limited balanced
truncation  methods in terms of absolute and relative errors for the two
applied input signals~\eqref{eqn:input_sco}, in Table~\ref{tab:sco_tl}, the
ROM sizes, absolute and relative errors are given.
One can observe that all ROMs are of order $4$, stable and
have locally significantly smaller errors than globally.

Again, the result of the two-step approaches are only marginal distinguishable
from the results of the full-order Gramians, where we used the global
sampling between $10^{-4}$ and $10^{4}$\,Hz to pre-approximate the system's
behavior.
Therefore, those results are also not shown here.

\begin{sidewaystable}[t]
  \renewcommand{\arraystretch}{1.25}
  \caption{Frequency-limited ROMs for the single chain oscillator
    (full-order Gramians).}
  \label{tab:sco_fl}
  \begin{center}
    \begin{tabular}{r|c|c|c|c|c|c|c|c}
      & \multicolumn{1}{c|}{\myP}
        & \multicolumn{1}{c|}{\myPM}
        & \multicolumn{1}{c|}{\myPV}
        & \multicolumn{1}{c|}{\myVP}
        & \multicolumn{1}{c|}{\myVPM}
        & \multicolumn{1}{c|}{\myV}
        & \multicolumn{1}{c|}{\myFV}
        & \multicolumn{1}{c}{\mySO}\\ \hline
      ROM sizes & 2 & 2 & 2 & 2 & 2 & 2 & 2 & 2 \\
      Stability & \cmark & \cmark & \cmark & \cmark & \cmark & \cmark
        & \cmark & \cmark \\ \hline
      Global absolute errors & \tabnum{1.011}{-}{01} & \tabnum{1.011}{-}{01}
        & \tabnum{1.011}{-}{01} & \tabnum{1.011}{-}{01}
        & \tabnum{1.011}{-}{01} & \tabnum{1.011}{-}{01}
        & \tabnum{1.012}{-}{01} & \tabnum{1.012}{-}{01} \\
      Local absolute errors & \tabnum{4.276}{-}{11} & \tabnum{4.277}{-}{11}
        & \tabnum{4.276}{-}{11} & \tabnum{7.439}{-}{11}
        & \tabnum{7.439}{-}{11} & \tabnum{7.439}{-}{11}
        & \tabnum{4.276}{-}{11} & \tabnum{7.439}{-}{11} \\ \hline
      Global relative errors & \tabnum{2.888}{-}{01} & \tabnum{2.888}{-}{01}
        & \tabnum{2.888}{-}{01} & \tabnum{2.888}{-}{01}
        & \tabnum{2.888}{-}{01} & \tabnum{2.888}{-}{01}
        & \tabnum{2.889}{-}{01} & \tabnum{2.889}{-}{01} \\
      Local relative errors & \tabnum{1.766}{-}{07} & \tabnum{1.766}{-}{07}
        & \tabnum{1.766}{-}{07} & \tabnum{3.072}{-}{07}
        & \tabnum{3.072}{-}{07} & \tabnum{3.072}{-}{07}
        & \tabnum{1.766}{-}{07} & \tabnum{3.072}{-}{07}
    \end{tabular}
  \end{center}
  
  \caption{Time-limited ROMs for the single chain oscillator
    (full-order Gramians).}
  \label{tab:sco_tl}
  \begin{center}
    \begin{tabular}{c|r|c|c|c|c|c|c|c|c}
      \multicolumn{2}{c}{} & \multicolumn{1}{|c|}{\myP}
        & \multicolumn{1}{c|}{\myPM}
        & \multicolumn{1}{c|}{\myPV}
        & \multicolumn{1}{c|}{\myVP}
        & \multicolumn{1}{c|}{\myVPM}
        & \multicolumn{1}{c|}{\myV}
        & \multicolumn{1}{c|}{\myFV}
        & \multicolumn{1}{c}{\mySO}\\ \hline
      \multicolumn{2}{r|}{ROM sizes} & 4 & 4 & 4 & 4 & 4 & 4 & 4 & 4 \\
      \multicolumn{2}{r|}{Stability} & \cmark & \cmark & \cmark & \cmark
        & \cmark & \cmark & \cmark & \cmark \\ \hline
      \multirow{4}{*}{$u_{\step}$} &
      Global absolute errors & \tabnum{9.621}{-}{04} & \tabnum{1.020}{-}{03}
        & \tabnum{9.619}{-}{04} & \tabnum{9.401}{-}{04}
        & \tabnum{9.985}{-}{04} & \tabnum{9.393}{-}{04}
        & \tabnum{9.880}{-}{04} & \tabnum{9.568}{-}{04} \\
      & Local absolute errors & \tabnum{7.953}{-}{07} & \tabnum{6.408}{-}{07}
        & \tabnum{7.980}{-}{07} & \tabnum{1.456}{-}{06}
        & \tabnum{2.866}{-}{06} & \tabnum{1.445}{-}{06}
        & \tabnum{8.597}{-}{07} & \tabnum{1.170}{-}{06} \\ \cline{2-10}
      & Global relative errors & \tabnum{4.724}{-}{03} & \tabnum{5.014}{-}{03}
        & \tabnum{4.723}{-}{03} & \tabnum{4.616}{-}{03}
        & \tabnum{4.910}{-}{03} & \tabnum{4.611}{-}{03}
        & \tabnum{4.853}{-}{03} & \tabnum{4.697}{-}{03} \\
      & Local relative errors & \tabnum{4.204}{-}{05} & \tabnum{1.256}{-}{05}
        & \tabnum{4.217}{-}{05} & \tabnum{4.617}{-}{05}
        & \tabnum{1.384}{-}{05} & \tabnum{4.634}{-}{05}
        & \tabnum{4.953}{-}{05} & \tabnum{4.503}{-}{05} \\ \hline
      \multirow{4}{*}{$u_{\sin}$} &
      Global absolute errors & \tabnum{5.232}{-}{05} & \tabnum{5.215}{-}{05}
        & \tabnum{5.231}{-}{05} & \tabnum{5.081}{-}{05}
        & \tabnum{5.045}{-}{05} & \tabnum{5.079}{-}{05}
        & \tabnum{5.350}{-}{05} & \tabnum{5.208}{-}{05} \\
      & Local absolute errors & \tabnum{8.600}{-}{07} & \tabnum{4.580}{-}{07}
        & \tabnum{8.619}{-}{07} & \tabnum{9.591}{-}{07}
        & \tabnum{4.961}{-}{07} & \tabnum{9.638}{-}{07}
        & \tabnum{9.471}{-}{07} & \tabnum{9.263}{-}{07} \\ \cline{2-10}
      & Global relative errors & \tabnum{8.275}{-}{02} & \tabnum{8.066}{-}{02}
        & \tabnum{8.273}{-}{02} & \tabnum{8.030}{-}{02}
        & \tabnum{7.827}{-}{02} & \tabnum{8.026}{-}{02}
        & \tabnum{8.465}{-}{02} & \tabnum{8.231}{-}{02} \\
      & Local relative errors & \tabnum{3.053}{-}{04} & \tabnum{1.150}{-}{04}
        & \tabnum{3.062}{-}{04} & \tabnum{3.264}{-}{04}
        & \tabnum{1.261}{-}{04} & \tabnum{3.284}{-}{04}
        & \tabnum{3.526}{-}{04} & \tabnum{3.214}{-}{04}
    \end{tabular}
  \end{center}
\end{sidewaystable}


\subsection{Crankshaft}

\begin{figure}[tb]
  \begin{center}
    \includegraphics[scale = 1]{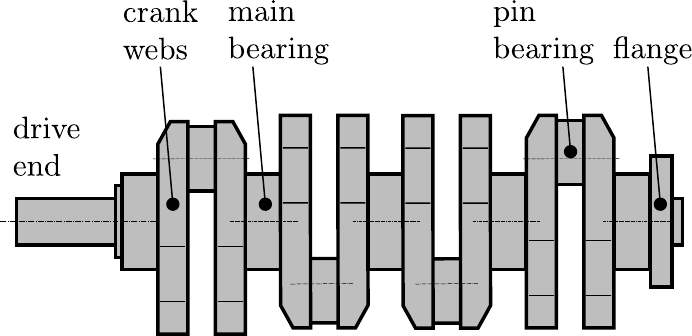}
  \end{center}
  \caption{Crankshaft of a four-cylinder engine~\cite{morNowKEetal13}.}
  \label{fig:cs}
\end{figure}

The crankshaft is a model from the University Stuttgart, describing the
crankshaft of a four-cylinder engine~\cite{morNowKEetal13}, which is shown
in Figure~\ref{fig:cs}.
After discretization by the finite element method, the constraint model is
of dimension $n = 42\,126$ with $m = p = 35$ inputs and outputs.
Due to the rigid elements, coupling the interface nodes, the system has
several eigenvalues at zero.
Therefore, we apply the shift $\alpha = 0.01$, as suggested in
Section~\ref{sec:alphashifts}, to make the system asymptotically stable during
the computations of the matrix equations and low-rank projection matrices.

\subsubsection{Frequency domain}

\begin{figure}[tb]
  \begin{center}
    \tikzexternalenable
    \tikzsetnextfilename{crankshaft_soflbt_a0_01}
    \input{graphics/crankshaft_soflbt_a0_01.tikz}
    \tikzexternaldisable
    \vspace{-\baselineskip}
  \end{center}
  \caption{Frequency-limited ROMs for the crankshaft
    (full-order Gramians).}
  \label{fig:cs_fl}
\end{figure}

\begin{figure}[tb]
  \begin{center}
    \tikzexternalenable
    \tikzsetnextfilename{crankshaft_soflbt_hybrid_local_a0_01}
    \input{graphics/crankshaft_soflbt_hybrid_local_a0_01.tikz}
    \tikzexternaldisable
    \vspace{-\baselineskip}
  \end{center}
  \caption{Frequency-limited ROMs for the crankshaft
    (two-step methods).}
  \label{fig:cs_flh}
\end{figure}
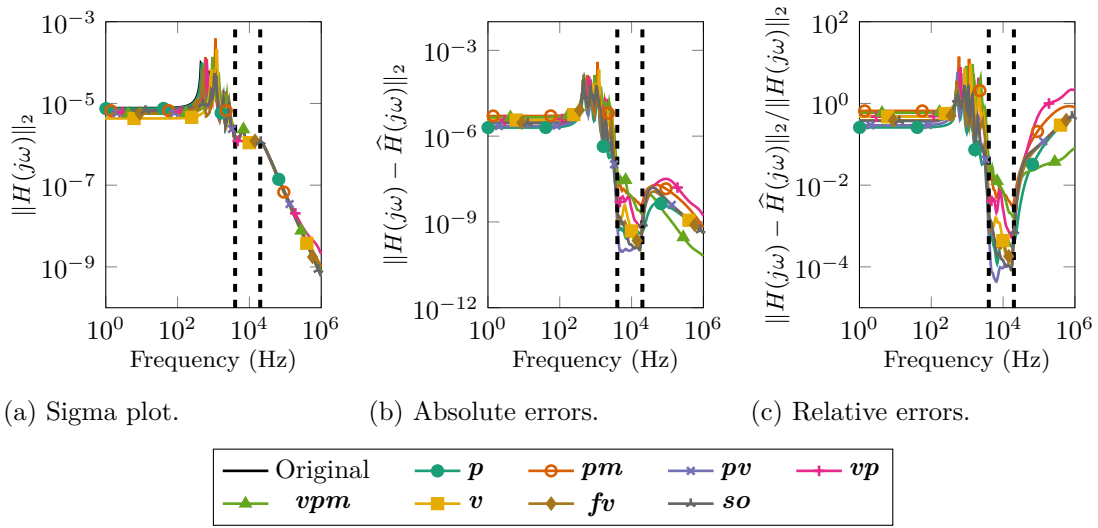

In the frequency domain, we are interested in the actual working range of the
crankshaft between $4$ and $20$\,kHz.
Figure~\ref{fig:cs_fl} shows the results for using the full-order
frequency-limited Gramians.
The frequency range of interest lies again between the two vertical
dashed lines.
We can see that all ROMs approximate the frequency region of
interest better than the global region.
Also Table~\ref{tab:cs_fl} shows the desired approximation behavior in terms 
of the errors.
In this example, some of the computed ROMs are unstable as denoted by x-marks
in Table~\ref{tab:cs_fl}.
It should be noted that even for the same order some methods might produce
unstable models while others do not.

In this example, we also applied the two-step approach with $200$ frequency
sample points in the region of interest to generate the intermediate model of
order $447$.
Those results can be seen in Figure~\ref{fig:cs_flh}.
Table~\ref{tab:cs_flh} shows that the ROMs produced by the two-step approach
are slightly larger in dimension and also partially in errors, while the same
methods (\myPM, \myVP, \myVPM, \mySO) as for the full-order Gramian approach
produce unstable models.

\subsubsection{Time domain}

\begin{figure}[tb]
  \begin{center}
    \tikzexternalenable
    \tikzsetnextfilename{crankshaft_sotlbt_mixed_a0_01}
    \input{graphics/crankshaft_sotlbt_mixed_a0_01.tikz}
    \tikzexternaldisable
    \vspace{-\baselineskip}
  \end{center}
  \caption{Absolute and relative errors of time-limited ROMs for the crankshaft
    with inputs $u_{\step}$ (a), (b) and $u_{\sin}$ (c), (d) (full-order
    Gramians).}
  \label{fig:cs_tl}
\end{figure}
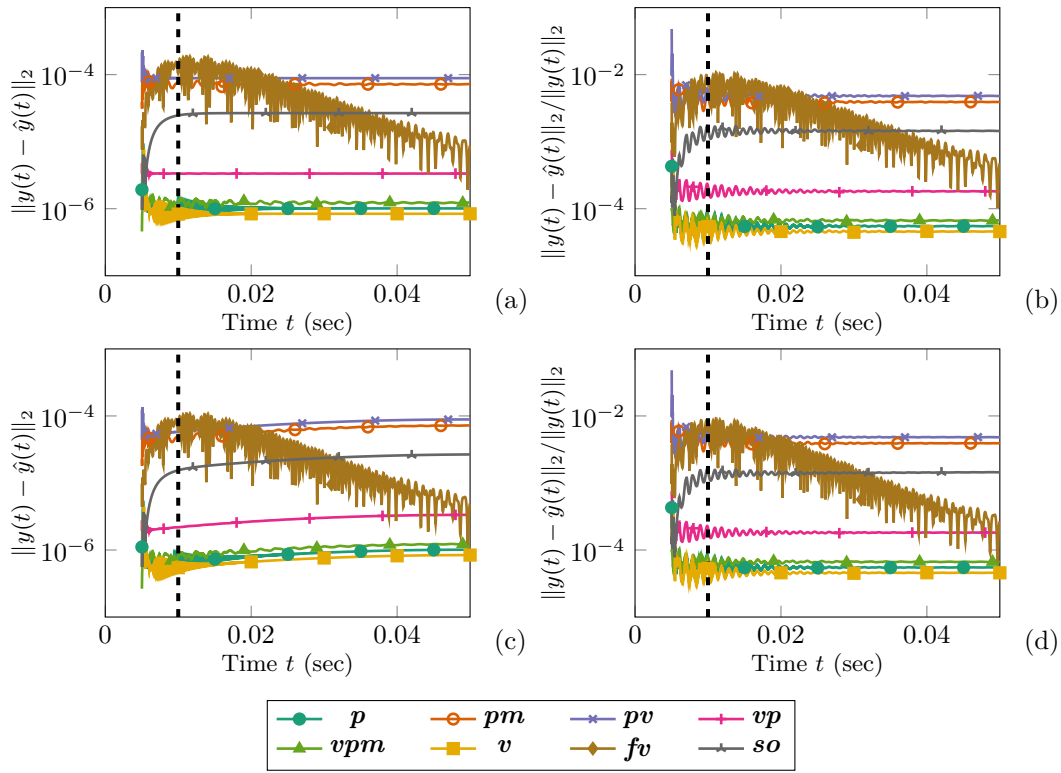

In the time domain, we consider just the first $0.01$\,s of using the
crankshaft, while the full simulation runs over a time range of $[0, 0.05]$\,s.
As test input signals, we apply
\begin{align*}
  \begin{aligned}
    u_{\step}(t) & = 3000 \delta(t - 0.005) \cdot \mathds{1}_{35} &&
      \text{and} & u_{\sin}(t) & = 1500(sin(10 \pi t) + 1)\delta(t - 0.005)
      \cdot \mathds{1}_{35},
  \end{aligned}
\end{align*}
where $\mathds{1}_{35}$ denotes the ones vector of length $35$.
The results for the time-limited balanced truncation with the full-order
Gramians can be seen in Figure~\ref{fig:cs_tl} and Table~\ref{tab:cs_tl}.
Only one unstable model (\myVPM) was computed, which still gives suitable
approximation results, and all ROMs have small enough errors in the time domain.
Even so, we recognize that the local approximation error is only in some
cases a bit smaller than the global one.

For the two-step approach, we computed $200$ logarithmically equidistant
distributed samples in the frequency domain between $10^{-2}$ and $10^{6}$\,Hz.
The intermediate model had the order $876$.
Since the resulting ROMs are of the same order as the ones computed via the 
full-order Gramians, featuring the same stability properties, and are only
slightly worse in terms of the time domain errors than in Table~\ref{tab:cs_tl},
we skip the additional presentation of those results here.

\begin{sidewaystable}[t]
  \renewcommand{\arraystretch}{1.25}
  \caption{Frequency-limited ROMs for the crankshaft
    (full-order Gramians).}
  \label{tab:cs_fl}
  \begin{center}
    \begin{tabular}{r|c|c|c|c|c|c|c|c}
      & \multicolumn{1}{c|}{\myP}
        & \multicolumn{1}{c|}{\myPM}
        & \multicolumn{1}{c|}{\myPV}
        & \multicolumn{1}{c|}{\myVP}
        & \multicolumn{1}{c|}{\myVPM}
        & \multicolumn{1}{c|}{\myV}
        & \multicolumn{1}{c|}{\myFV}
        & \multicolumn{1}{c}{\mySO}\\ \hline
      ROM sizes & 77 & 77 & 65 & 88 & 88 & 69 & 77 & 77 \\
      Stability & \cmark & \xmark & \cmark & \xmark & \xmark & \cmark
        & \cmark & \xmark \\ \hline
      Global absolute errors & \tabnum{9.367}{-}{05} & \tabnum{3.237}{-}{04}
        & \tabnum{9.280}{-}{05} & \tabnum{1.141}{-}{04}
        & \tabnum{9.601}{-}{05} & \tabnum{9.265}{-}{05}
        & \tabnum{9.361}{-}{05} & \tabnum{9.345}{-}{05} \\
      Local absolute errors & \tabnum{1.588}{-}{09} & \tabnum{1.816}{-}{08}
        & \tabnum{9.855}{-}{10} & \tabnum{5.497}{-}{09}
        & \tabnum{4.978}{-}{08} & \tabnum{4.413}{-}{10}
        & \tabnum{1.011}{-}{10} & \tabnum{2.818}{-}{10} \\ \hline
      Global relative errors & \tabnum{4.627}{+}{00} & \tabnum{2.082}{+}{01}
        & \tabnum{2.353}{+}{00} & \tabnum{1.439}{+}{01}
        & \tabnum{4.682}{+}{00} & \tabnum{4.718}{+}{00}
        & \tabnum{3.963}{+}{00} & \tabnum{2.652}{+}{00} \\
      Local relative errors & \tabnum{1.327}{-}{03} & \tabnum{1.754}{-}{02}
        & \tabnum{8.237}{-}{04} & \tabnum{5.117}{-}{03}
        & \tabnum{4.807}{-}{02} & \tabnum{4.261}{-}{04}
        & \tabnum{9.759}{-}{05} & \tabnum{2.722}{-}{04}
    \end{tabular}
  \end{center}
  
  \caption{Frequency-limited ROMs for the crankshaft (two-step methods).}
  \label{tab:cs_flh}
  \begin{center}
    \begin{tabular}{r|c|c|c|c|c|c|c|c}
      & \multicolumn{1}{c|}{\myP}
        & \multicolumn{1}{c|}{\myPM}
        & \multicolumn{1}{c|}{\myPV}
        & \multicolumn{1}{c|}{\myVP}
        & \multicolumn{1}{c|}{\myVPM}
        & \multicolumn{1}{c|}{\myV}
        & \multicolumn{1}{c|}{\myFV}
        & \multicolumn{1}{c}{\mySO}\\ \hline
      ROM sizes & 84 & 84 & 67 & 93 & 93 & 70 & 84 & 70 \\
      Stability & \cmark & \xmark & \cmark & \xmark & \xmark & \cmark
        & \cmark & \xmark \\ \hline
      Global absolute errors & \tabnum{1.405}{-}{04} & \tabnum{3.945}{-}{04}
        & \tabnum{1.026}{-}{04} & \tabnum{1.204}{-}{04}
        & \tabnum{1.364}{-}{04} & \tabnum{2.138}{-}{04}
        & \tabnum{1.065}{-}{04} & \tabnum{9.057}{-}{05} \\
      Local absolute errors & \tabnum{2.037}{-}{09} & \tabnum{2.569}{-}{08}
        & \tabnum{8.225}{-}{10} & \tabnum{7.911}{-}{09}
        & \tabnum{3.709}{-}{08} & \tabnum{3.846}{-}{09}
        & \tabnum{1.297}{-}{09} & \tabnum{1.774}{-}{09} \\ \hline
      Global relative errors & \tabnum{2.041}{+}{00} & \tabnum{1.400}{+}{01}
        & \tabnum{5.712}{+}{00} & \tabnum{7.743}{+}{00}
        & \tabnum{9.187}{+}{00} & \tabnum{9.393}{+}{00}
        & \tabnum{3.565}{+}{00} & \tabnum{2.865}{+}{00} \\
      Local relative errors & \tabnum{1.967}{-}{03} & \tabnum{2.481}{-}{02}
        & \tabnum{6.874}{-}{04} & \tabnum{7.775}{-}{03}
        & \tabnum{3.462}{-}{02} & \tabnum{1.810}{-}{03}
        & \tabnum{1.252}{-}{03} & \tabnum{1.713}{-}{03}
    \end{tabular}
  \end{center}
\end{sidewaystable}


\subsection{Artificial fishtail}

\begin{figure}[tb]
  \begin{center}
    \includegraphics[scale = .3]{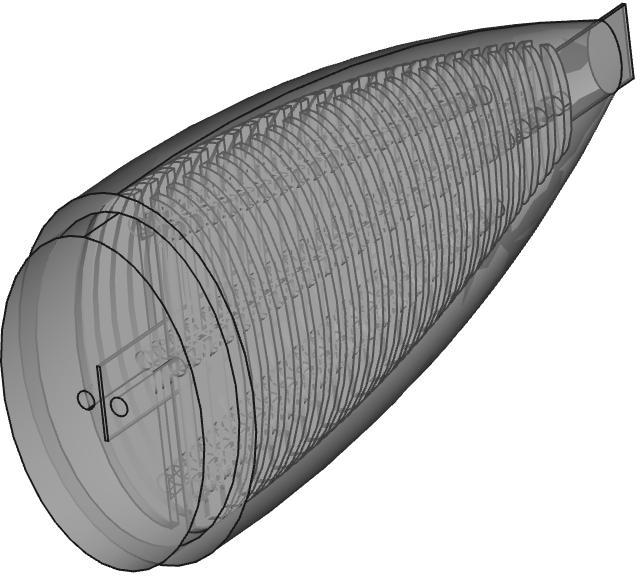}
  \end{center}
  \caption{Transparent sketch of the artificial fishtail with embedded fluid
    chambers.}
  \label{fig:ft}
\end{figure}

The artificial fishtail is a mechanical system, describing the movement of a
fishtail-shaped structure by using the fluid elastomer actuation principle.
Figure~\ref{fig:ft} shows a transparent sketch of the fishtail model consisting
of a carbon beam in the center and a silicon hull around.
A more detailed description of the model as well as a comparison of
structure-preserving second-order model reduction techniques for this example
can be found in~\cite{morSaaSW19}.
After spatial discretization by the finite element method, the resulting
second-order system has $n = 779\,232$ states describing the model.
By the actuation principle, we have $m = 1$ input and a sensor is measuring the
displacement of the fishtail's tip in all spatial dimensions, i.e., we have
$p = 3$ position outputs and no velocity outputs.
The discretized data is available as open benchmark at~\cite{SieKM19}.
The computations were done without an $\alpha$-shift ($\alpha = 0$).

\subsubsection{Frequency domain}

\begin{figure}[tb]
  \begin{center}
    \tikzexternalenable
    \tikzsetnextfilename{fish_tail_soflbt_a0}
    \input{graphics/fish_tail_soflbt_a0.tikz}
    \tikzexternaldisable
    \vspace{-2\baselineskip}
  \end{center}
  \caption{Frequency-limited ROMs for the artificial fishtail
    (full-order Gramians).}
  \label{fig:ft_fl}
\end{figure}
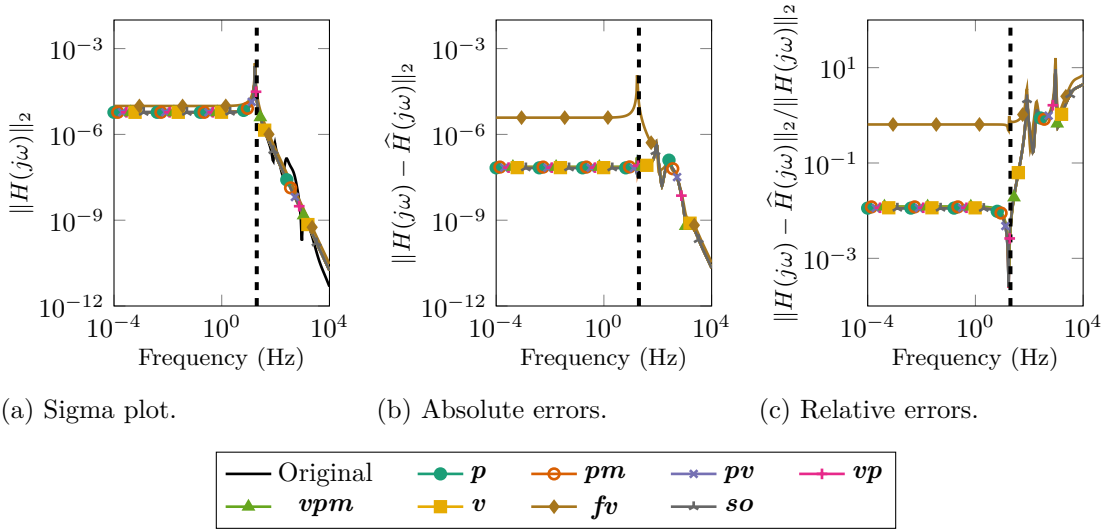

In the frequency domain, the range of interest for the fishtail model lies
between $0$ and $20$\,Hz, since higher frequencies are physically not realizable.
Figure~\ref{fig:ft_fl} shows the results for the frequency-limited balanced
truncation methods, based on the full-order Gramians.
Except for the \myFV{} balancing there is no visible difference between
the ROMs and the full-order model.
The error plots show that the approximation reached a sufficiently small error
in the region of interest.
Table~\ref{tab:ft_fl} shows the corresponding maximum absolute and relative
error in the local and global frequency regions.
It is remarkable that the methods were able to approximate the original model,
having around $780\,000$ states, by stable order $1$ systems in the region of
interest.
While the absolute errors are comparable between local and global region,
the relative errors show again the strength of the frequency-limited method.

\subsubsection{Time domain}

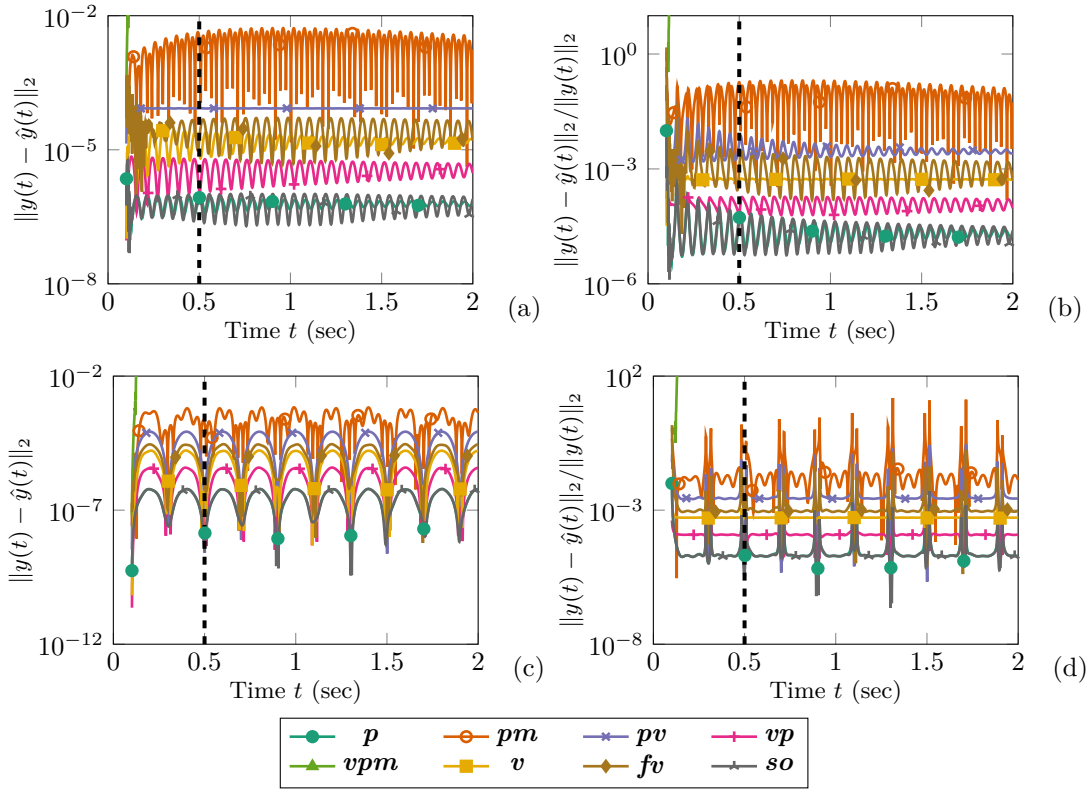
\begin{figure}[tb]
  \begin{center}
    \tikzexternalenable
    \tikzsetnextfilename{fish_tail_sotlbt_a0}
    \input{graphics/fish_tail_sotlbt_a0.tikz}
    \tikzexternaldisable
    \vspace{-\baselineskip}
  \end{center}
  \caption{Absolute and relative errors of time-limited ROMs for the artificial
    fishtail with inputs $u_{\step}$ (a), (b) and $u_{\sin}$ (c), (d)
    (full-order Gramians).}
  \label{fig:ft_tl}
\end{figure}

In the time domain, the fishtail is simulated from $0$ to $2$\,s.
For our time-limited methods we consider the time range up to $0.5$\,s and
as inputs, the following two signals are considered
\begin{align*}
  \begin{aligned}
    u_{\step}(t) & = 5000 \delta(t - 0.1) && \text{and} &
      u_{\sin}(t) & = 2500(sin(10 \pi (t - 1.35)) + 1)\delta(t - 0.1).
  \end{aligned}
\end{align*}
Figure~\ref{fig:ft_tl} and Table~\ref{tab:ft_tl} show the results.
Except for the models generated by \myPM, \myVPM{} and \myFV, the computed ROMs
have acceptable small errors in the time domain.
Also, only the \myVPM{} ROM is unstable.
The errors in the local region are sometimes a bit smaller than the global one 
as we were aiming for by the method.

The two-step approach here used $200$ logarithmically equidistant sample points
in the frequency range from $10^{-4}$ to $10^{4}$\,Hz, which gave an
intermediate model of order $100$.
The results of the ROMs computed by the two-step approach differ a bit from the
ones generated by the full-order Gramians.
Those results can be seen in Table~\ref{tab:ft_tlh}.
There, shown errors are partially smaller or larger than in Table~\ref{tab:ft_tl}
and also we note that for the two-step approach, the \myVPM{} model is also
unstable but still gives usable results for both applied input signals.

\begin{sidewaystable}[t]
  \renewcommand{\arraystretch}{1.25}
  \caption{Time-limited ROMs for the crankshaft (full-order Gramians).}
  \label{tab:cs_tl}
  \begin{center}
    \begin{tabular}{c|r|c|c|c|c|c|c|c|c}
      \multicolumn{2}{c}{} & \multicolumn{1}{|c|}{\myP}
        & \multicolumn{1}{c|}{\myPM}
        & \multicolumn{1}{c|}{\myPV}
        & \multicolumn{1}{c|}{\myVP}
        & \multicolumn{1}{c|}{\myVPM}
        & \multicolumn{1}{c|}{\myV}
        & \multicolumn{1}{c|}{\myFV}
        & \multicolumn{1}{c}{\mySO}\\ \hline
      \multicolumn{2}{r|}{ROM sizes} & 58 & 58 & 37 & 132 & 132 & 59 & 58
        & 59 \\
      \multicolumn{2}{r|}{Stability} & \cmark & \cmark & \cmark & \cmark
        & \xmark & \cmark & \cmark & \cmark \\ \hline
      \multirow{4}{*}{$u_{\step}$} &
      Global absolute errors & \tabnum{9.442}{-}{06} & \tabnum{8.765}{-}{05}
        & \tabnum{2.331}{-}{04} & \tabnum{4.771}{-}{06}
        & \tabnum{1.567}{-}{06} & \tabnum{7.733}{-}{06}
        & \tabnum{1.884}{-}{04} & \tabnum{2.669}{-}{05} \\
      & Local absolute errors & \tabnum{9.442}{-}{06} & \tabnum{8.765}{-}{05}
        & \tabnum{2.331}{-}{04} & \tabnum{4.771}{-}{06}
        & \tabnum{1.567}{-}{06} & \tabnum{7.733}{-}{06}
        & \tabnum{1.554}{-}{04} & \tabnum{2.467}{-}{05} \\ \cline{2-10}
      & Global relative errors & \tabnum{5.573}{-}{04} & \tabnum{8.707}{-}{03}
        & \tabnum{4.803}{-}{02} & \tabnum{7.460}{-}{04}
        & \tabnum{1.039}{-}{04} & \tabnum{4.103}{-}{04}
        & \tabnum{1.078}{-}{02} & \tabnum{1.731}{-}{03} \\
      & Local relative errors & \tabnum{5.573}{-}{04} & \tabnum{8.707}{-}{03}
        & \tabnum{4.803}{-}{02} & \tabnum{7.460}{-}{04}
        & \tabnum{1.039}{-}{04} & \tabnum{4.103}{-}{04}
        & \tabnum{9.845}{-}{03} & \tabnum{1.731}{-}{03} \\ \hline
      \multirow{4}{*}{$u_{\sin}$} &
      Global absolute errors & \tabnum{5.459}{-}{06} & \tabnum{7.231}{-}{05}
        & \tabnum{1.349}{-}{04} & \tabnum{3.345}{-}{06}
        & \tabnum{1.233}{-}{06} & \tabnum{4.472}{-}{06}
        & \tabnum{1.089}{-}{04} & \tabnum{2.664}{-}{05} \\
      & Local absolute errors & \tabnum{5.459}{-}{06} & \tabnum{5.245}{-}{05}
        & \tabnum{1.349}{-}{04} & \tabnum{2.760}{-}{06}
        & \tabnum{9.133}{-}{07} & \tabnum{4.472}{-}{06}
        & \tabnum{8.980}{-}{05} & \tabnum{1.547}{-}{05} \\ \cline{2-10}
      & Global relative errors & \tabnum{5.559}{-}{04} & \tabnum{8.707}{-}{03}
        & \tabnum{4.803}{-}{02} & \tabnum{7.460}{-}{04}
        & \tabnum{1.028}{-}{04} & \tabnum{4.103}{-}{04}
        & \tabnum{9.189}{-}{03} & \tabnum{1.610}{-}{03} \\
      & Local relative errors & \tabnum{5.559}{-}{04} & \tabnum{8.707}{-}{03}
        & \tabnum{4.803}{-}{02} & \tabnum{7.460}{-}{04}
        & \tabnum{1.028}{-}{04} & \tabnum{4.103}{-}{04}
        & \tabnum{8.605}{-}{03} & \tabnum{1.610}{-}{03}
    \end{tabular}
  \end{center}
  
  \caption{Frequency-limited ROMs for the artificial fishtail
    (full-order Gramians).}
  \label{tab:ft_fl}
  \begin{center}
    \begin{tabular}{r|c|c|c|c|c|c|c|c}
      & \multicolumn{1}{|c|}{\myP}
        & \multicolumn{1}{c|}{\myPM}
        & \multicolumn{1}{c|}{\myPV}
        & \multicolumn{1}{c|}{\myVP}
        & \multicolumn{1}{c|}{\myVPM}
        & \multicolumn{1}{c|}{\myV}
        & \multicolumn{1}{c|}{\myFV}
        & \multicolumn{1}{c}{\mySO}\\ \hline
      ROM sizes & 1 & 1 & 1 & 1 & 1 & 1 & 1 & 1 \\
      Stability & \cmark & \cmark & \cmark & \cmark & \cmark & \cmark
        & \cmark & \cmark \\ \hline
      Global absolute errors & \tabnum{4.409}{-}{07} & \tabnum{4.409}{-}{07}
        & \tabnum{4.409}{-}{07} & \tabnum{4.409}{-}{07}
        & \tabnum{4.409}{-}{07} & \tabnum{4.409}{-}{07}
        & \tabnum{1.172}{-}{04} & \tabnum{4.409}{-}{07} \\
      Local absolute errors & \tabnum{1.046}{-}{07} & \tabnum{1.538}{-}{07}
        & \tabnum{1.043}{-}{07} & \tabnum{8.975}{-}{08}
        & \tabnum{1.558}{-}{07} & \tabnum{8.964}{-}{08}
        & \tabnum{1.172}{-}{04} & \tabnum{1.045}{-}{07} \\ \hline
      Global relative errors & \tabnum{9.182}{+}{00} & \tabnum{9.176}{+}{00}
        & \tabnum{9.182}{+}{00} & \tabnum{9.181}{+}{00}
        & \tabnum{9.174}{+}{00} & \tabnum{9.181}{+}{00}
        & \tabnum{1.596}{+}{01} & \tabnum{9.182}{+}{00} \\
      Local relative errors & \tabnum{1.132}{-}{02} & \tabnum{1.200}{-}{02}
        & \tabnum{1.132}{-}{02} & \tabnum{1.150}{-}{02}
        & \tabnum{1.219}{-}{02} & \tabnum{1.150}{-}{02}
        & \tabnum{9.557}{-}{01} & \tabnum{1.132}{-}{02}
    \end{tabular}
  \end{center}
\end{sidewaystable}
  
\begin{sidewaystable}[t]
  \renewcommand{\arraystretch}{1.25}
  \caption{Time-limited ROMs for the artificial fishtail (full-order Gramians).}
  \label{tab:ft_tl}
  \begin{center}
    \begin{tabular}{c|r|c|c|c|c|c|c|c|c}
      \multicolumn{2}{c}{} & \multicolumn{1}{|c|}{\myP}
        & \multicolumn{1}{c|}{\myPM}
        & \multicolumn{1}{c|}{\myPV}
        & \multicolumn{1}{c|}{\myVP}
        & \multicolumn{1}{c|}{\myVPM}
        & \multicolumn{1}{c|}{\myV}
        & \multicolumn{1}{c|}{\myFV}
        & \multicolumn{1}{c}{\mySO}\\ \hline
      \multicolumn{2}{r|}{ROM sizes} & 4 & 4 & 2 & 6 & 6 & 4 & 4 & 4 \\
      \multicolumn{2}{r|}{Stability} & \cmark & \cmark & \cmark & \cmark
        & \xmark & \cmark & \cmark & \cmark \\ \hline
      \multirow{4}{*}{$u_{\step}$} &
      Global absolute errors & \tabnum{5.523}{-}{06} & \tabnum{5.277}{-}{03}
        & \tabnum{2.320}{-}{04} & \tabnum{7.032}{-}{06}
        & $\infty$ & \tabnum{3.049}{-}{05}
        & \tabnum{4.650}{-}{04} & \tabnum{6.087}{-}{06} \\
      & Local absolute errors & \tabnum{5.523}{-}{06} & \tabnum{4.282}{-}{03}
        & \tabnum{2.320}{-}{04} & \tabnum{7.032}{-}{06}
        & $\infty$ & \tabnum{3.049}{-}{05}
        & \tabnum{4.650}{-}{04} & \tabnum{6.087}{-}{06} \\ \cline{2-10}
      & Global relative errors & \tabnum{9.961}{-}{03} & \tabnum{4.577}{-}{01}
        & \tabnum{1.524}{-}{01} & \tabnum{4.127}{-}{04}
        & $\infty$ & \tabnum{2.799}{-}{03}
        & \tabnum{1.489}{+}{00} & \tabnum{8.162}{-}{03} \\
      & Local relative errors & \tabnum{9.961}{-}{03} & \tabnum{4.577}{-}{01}
        & \tabnum{1.524}{-}{01} & \tabnum{4.127}{-}{04}
        & $\infty$ & \tabnum{2.799}{-}{03}
        & \tabnum{1.489}{+}{00} & \tabnum{8.162}{-}{03} \\ \hline
      \multirow{4}{*}{$u_{\sin}$} &
      Global absolute errors & \tabnum{6.103}{-}{07} & \tabnum{6.845}{-}{04}
        & \tabnum{8.434}{-}{05} & \tabnum{3.878}{-}{06}
        & $\infty$ & \tabnum{1.681}{-}{05}
        & \tabnum{2.898}{-}{05} & \tabnum{6.237}{-}{07} \\
      & Local absolute errors & \tabnum{6.094}{-}{07} & \tabnum{6.278}{-}{04}
        & \tabnum{8.434}{-}{05} & \tabnum{3.850}{-}{06}
        & $\infty$ & \tabnum{1.681}{-}{05}
        & \tabnum{2.846}{-}{05} & \tabnum{6.129}{-}{07} \\ \cline{2-10}
      & Global relative errors & \tabnum{9.961}{-}{03} & \tabnum{1.525}{+}{01}
        & \tabnum{2.819}{-}{01} & \tabnum{6.047}{-}{03}
        & $\infty$ & \tabnum{3.089}{-}{03}
        & \tabnum{1.489}{+}{00} & \tabnum{8.162}{-}{03} \\
      & Local relative errors & \tabnum{9.961}{-}{03} & \tabnum{1.549}{+}{00}
        & \tabnum{1.224}{-}{01} & \tabnum{8.350}{-}{04}
        & $\infty$ & \tabnum{1.192}{-}{03}
        & \tabnum{1.489}{+}{00} & \tabnum{8.162}{-}{03}
    \end{tabular}
  \end{center}
  
  \caption{Time-limited ROMs for the artificial fishtail (two-step methods).}
  \label{tab:ft_tlh}
  \begin{center}
    \begin{tabular}{c|r|c|c|c|c|c|c|c|c}
      \multicolumn{2}{c}{} & \multicolumn{1}{|c|}{\myP}
        & \multicolumn{1}{c|}{\myPM}
        & \multicolumn{1}{c|}{\myPV}
        & \multicolumn{1}{c|}{\myVP}
        & \multicolumn{1}{c|}{\myVPM}
        & \multicolumn{1}{c|}{\myV}
        & \multicolumn{1}{c|}{\myFV}
        & \multicolumn{1}{c}{\mySO}\\ \hline
      \multicolumn{2}{r|}{ROM sizes} & 4 & 4 & 2 & 9 & 9 & 4 & 4 & 4 \\
      \multicolumn{2}{r|}{Stability} & \cmark & \cmark & \cmark & \cmark
        & \xmark & \cmark & \cmark & \cmark \\ \hline
      \multirow{4}{*}{$u_{\step}$} &
      Global absolute errors & \tabnum{5.506}{-}{06} & \tabnum{1.306}{-}{03}
        & \tabnum{2.308}{-}{04} & \tabnum{6.210}{-}{06}
        & \tabnum{1.394}{-}{03} & \tabnum{6.649}{-}{05}
        & \tabnum{2.229}{-}{04} & \tabnum{7.206}{-}{06} \\
      & Local absolute errors & \tabnum{5.506}{-}{06} & \tabnum{1.306}{-}{03}
        & \tabnum{2.308}{-}{04} & \tabnum{6.210}{-}{06}
        & \tabnum{1.137}{-}{03} & \tabnum{6.649}{-}{05}
        & \tabnum{2.229}{-}{04} & \tabnum{7.206}{-}{06} \\ \cline{2-10}
      & Global relative errors & \tabnum{1.088}{-}{02} & \tabnum{2.335}{+}{00}
        & \tabnum{1.517}{-}{01} & \tabnum{2.461}{-}{03}
        & \tabnum{2.321}{-}{01} & \tabnum{2.253}{-}{02}
        & \tabnum{9.866}{-}{01} & \tabnum{4.656}{-}{03} \\
      & Local relative errors & \tabnum{1.088}{-}{02} & \tabnum{2.335}{+}{00}
        & \tabnum{1.517}{-}{01} & \tabnum{2.461}{-}{03}
        & \tabnum{2.321}{-}{01} & \tabnum{2.253}{-}{02}
        & \tabnum{9.866}{-}{01} & \tabnum{4.656}{-}{03} \\ \hline
      \multirow{4}{*}{$u_{\sin}$} &
      Global absolute errors & \tabnum{9.836}{-}{07} & \tabnum{9.156}{-}{04}
        & \tabnum{8.371}{-}{05} & \tabnum{5.547}{-}{07}
        & \tabnum{3.619}{-}{04} & \tabnum{2.309}{-}{05}
        & \tabnum{9.389}{-}{07} & \tabnum{9.808}{-}{07} \\
      & Local absolute errors & \tabnum{9.885}{-}{07} & \tabnum{9.156}{-}{04}
        & \tabnum{8.371}{-}{05} & \tabnum{5.560}{-}{07}
        & \tabnum{3.887}{-}{04} & \tabnum{2.316}{-}{05}
        & \tabnum{9.389}{-}{07} & \tabnum{9.899}{-}{07} \\ \cline{2-10}
      & Global relative errors & \tabnum{1.088}{-}{02} & \tabnum{2.335}{+}{00}
        & \tabnum{2.775}{-}{01} & \tabnum{2.613}{-}{03}
        & \tabnum{4.312}{+}{00} & \tabnum{1.152}{-}{01}
        & \tabnum{9.866}{-}{01} & \tabnum{4.656}{-}{03} \\
      & Local relative errors & \tabnum{1.088}{-}{02} & \tabnum{2.335}{+}{00}
        & \tabnum{1.218}{-}{01} & \tabnum{1.585}{-}{03}
        & \tabnum{6.518}{-}{01} & \tabnum{3.076}{-}{02}
        & \tabnum{9.866}{-}{01} & \tabnum{4.656}{-}{03}
    \end{tabular}
  \end{center}
\end{sidewaystable}


\section{Conclusions}
\label{sec:conclusions}

We extended the frequency- and time-limited balanced truncation methods from
first-order systems to the second-order case by applying the different
second-order balancing approaches from the literature.
For the application of the introduced theory, we investigated
numerical methods for approximating the solution of the arising large-scale
sparse matrix equations with function right hand-sides as well as techniques
to deal with the difficulties arising from the second-order system structure.
The numerical examples show that the methods work for the purpose of
limited model reduction in the frequency domain and also for some examples in
time domain.
By comparison of the different balancing formulas, it was not possible to
determine a clear winner or loser.
Depending on the example, different balancing techniques performed better or
worse than the others.
Also, stability preservation is still an open problem for this type of model
reduction techniques, where we pointed out that the known modifications from
the first-order case are not necessarily stability preserving for second-order
systems.


\addcontentsline{toc}{section}{Acknowledgment}
\section*{Acknowledgment}
This work was supported by the German Research Foundation (DFG) Research
Training Group 2297 ``MathCoRe'', Magdeburg, and the German Research Foundation
(DFG) Priority Program 1897: ``Calm, Smooth and Smart -- Novel Approaches for
Influencing Vibrations by Means of Deliberately Introduced Dissipation''.

We would like to thank Patrick K{\"u}rschner who helped with an initial
version of the codes in the limited balanced truncation for large-scale sparse
second-order systems package~\cite{morBenW20}.


\addcontentsline{toc}{section}{References}
\bibliographystyle{plainurl}
\bibliography{bibtex/myref}

\end{document}

%% file: graphics/damped_mass_spring.tikz
%
%

\begin{tikzpicture}[
  every node/.style = {
    draw, 
    outer sep = 0pt,
    thick}
]

  \tikzstyle{damper} = [
    thick,
    decoration = {
      markings,  
      mark connection node = dmp,
      mark = at position 0.5 with {
        \node (dmp) [thick, inner sep = 0pt, transform shape, rotate = -90, 
          minimum width = 15pt, minimum height = 3pt, draw = none] {};
        \draw [thick] ($(dmp.north east)+(2pt,0)$) -- (dmp.south east) -- 
          (dmp.south west) -- ($(dmp.north west)+(2pt,0)$);
        \draw [thick] ($(dmp.north)+(0,-5pt)$) -- ($(dmp.north)+(0,5pt)$);
      }}, 
    decorate
  ]
  
  \tikzstyle{ground} = [
    fill,
    pattern = north east lines,
    draw = none,
    minimum width = 1cm,
    minimum height = 0.3cm,
    rotate = -90
  ]

  \tikzstyle{mass} = [
    circle,
    minimum size = 1cm
  ]
  
  \tikzstyle{spring} = [
    thick, 
    decorate, 
    decoration = {
      zigzag, 
      pre length = 0.2cm, 
      post length = 0.2cm, 
      segment length = 6,
      amplitude = 0.15cm}
  ]
  
  \newcommand*{\DrawSpringDamper}[4]{
    \draw[thick] #1 -- ++(0.2cm,0) coordinate (w);
    \draw[thick, line cap = rect] (w) -- ++(0,0.3cm) coordinate (ws);
    \draw[spring] (ws) -- ++(1.1cm,0) coordinate (sm) 
      node[midway, above = 0.12cm, draw = none] {#3};
    \draw[thick, line cap = rect] (sm) -- ++(0,-0.3cm) coordinate (m);
    \draw[thick, line cap = rect] (w) -- ++(0,-0.3cm) coordinate (wd);
    \draw[damper] (wd) -- ++(1.1cm,0) coordinate (dm) 
      node[midway, below = 0.22cm, draw = none] {#4};
    \draw[thick, line cap = rect] (dm) -- ++(0,0.3cm);
    \draw[thick] (m) -- #2;
  }
  
  \node[ground, fill = none] (wall0) {};
  
  \node[mass] (mass1) [right = 1.0cm of wall0.north] {$m_{1}$};
  \draw[-latex, ultra thick] (mass1.west) ++(-1cm,0) -- +(1cm,0) 
    node[midway, below, draw = none] {$u$};

  \node[ground] (wall1) 
    [below = 0.85cm of wall0.east, anchor = west] {};
  \draw[thick] (wall1.north west) -- (wall1.north east);
  \draw[thick, line cap = rect] (mass1.south) -- ++(0,-1.35cm) coordinate (me);
  \DrawSpringDamper
    {(wall1.north)}
    {(me)}
    {$\kappa_{1}$}
    {$\delta_{1}$}
  
  \DrawSpringDamper
    {(mass1.east)}
    {++(0.2cm,0) coordinate (lp)}
    {$k_{1}$}
    {$d_{1}$}
  \draw[thick, dash pattern = on \pgflinewidth off 0.1cm] (lp) -- ++(0.5cm,0)
    coordinate (lp);
  \DrawSpringDamper
    {(lp)}
    {++(0.2cm,0) coordinate (lp)}
    {$k_{i-1}$}
    {$d_{i-1}$}
  
  \node[mass, anchor = west] (massi) at (lp) {$m_{i}$};
  
  \node[ground] (walli) [right = 4.0cm of wall1.north, anchor = north] {};
  \draw[thick] (walli.north west) -- (walli.north east);
  \draw[thick, line cap = rect] (massi.south) -- ++(0,-1.35cm) coordinate (me);
  \DrawSpringDamper
    {(walli.north)}
    {(me)}
    {$\kappa_{i}$}
    {$\delta_{i}$}
  
  \DrawSpringDamper
    {(massi.east)}
    {++(0.2cm,0) coordinate (lp)}
    {$k_{i}$}
    {$d_{i}$}
  \draw[thick, dash pattern = on \pgflinewidth off 0.1cm] (lp) -- ++(0.5cm,0)
    coordinate (lp);
  \DrawSpringDamper
    {(lp)}
    {++(0.2cm,0) coordinate (lp)}
    {$k_{n-1}$}
    {$d_{n-1}$}
  
  \node[mass, anchor = west] (massg) at (lp) {$m_{n}$};
  
  \node[ground] (wallg) [right = 4.5cm of walli.north, anchor = north] {};
  \draw[thick] (wallg.north west) -- (wallg.north east);
  \draw[thick, line cap = rect] (massg.south) -- ++(0,-1.35cm) coordinate (me);
  \DrawSpringDamper
    {(wallg.north)}
    {(me)}
    {$\kappa_{n}$}
    {$\delta_{n}$}
  
  
\end{tikzpicture}

%% file: graphics/singlechain_soflbt_a0.tikz
\begin{tikzpicture}[every axis/.append style={cycle list name = freqlist}]
  \pgfplotstableread{graphics/data/singlechain_tf.dat}\tableFOM
  \pgfplotstableread{graphics/data/singlechain_soflbt_a0_tf.dat}\tableROM
  \pgfplotstableread{graphics/data/singlechain_soflbt_a0_abserr.dat}\tableABS
  \pgfplotstableread{graphics/data/singlechain_soflbt_a0_relerr.dat}\tableREL
  
  \node(tf){
    \begin{tikzpicture}
      \begin{loglogaxis}[%
        width  = .189\textwidth,
        height = .16\textheight,
        scale only axis,
        xmin = 1e-4,
        xmax = 1e+4,
        ymin = 1e-12,
        ymax = 1e+0,
        xminorticks = false,
        xlabel={\small Frequency (Hz)},
        xlabel style = {yshift = .3em},
        ylabel = {\small $\lVert H(j\omega) \rVert_{2}$},
        ylabel style = {yshift = -.5em}]
        
        \addplot table[x index=0, y index=1] {\tableFOM};
        \addplot table[x index=0, y index=1] {\tableROM};
        \addplot table[x index=0, y index=2] {\tableROM};
        \addplot table[x index=0, y index=3] {\tableROM};
        \addplot table[x index=0, y index=4] {\tableROM};
        \addplot table[x index=0, y index=5] {\tableROM};
        \addplot table[x index=0, y index=6] {\tableROM};
        \addplot table[x index=0, y index=7] {\tableROM};
        \addplot table[x index=0, y index=8] {\tableROM};
        
        \addplot[rngLine] coordinates {(1e0, 1e-12) (1e0, 1e0)};
        \addplot[rngLine] coordinates {(1e2, 1e-12) (1e2, 1e0)};
      \end{loglogaxis}
    \end{tikzpicture}};
    
  \node(abs)[right = 0cm of tf.south east, anchor = south west]{
    \begin{tikzpicture}
      \begin{loglogaxis}[%
        width  = .189\textwidth,
        height = .16\textheight,
        scale only axis,
        xmin = 1e-4,
        xmax = 1e+4,
        ymin = 1e-20,
        ymax = 1e+0,
        xminorticks = false,
        xlabel={\small Frequency (Hz)},
        xlabel style = {yshift = .3em},
        ylabel = {\small $\lVert H(j\omega) - \hH(j\omega) \rVert_{2}$},
        ylabel style = {yshift = -.5em}]
        
        \pgfplotsset{cycle list shift=1}
        \addplot table[x index=0, y index=1] {\tableABS};
        \addplot table[x index=0, y index=2] {\tableABS};
        \addplot table[x index=0, y index=3] {\tableABS};
        \addplot table[x index=0, y index=4] {\tableABS};
        \addplot table[x index=0, y index=5] {\tableABS};
        \addplot table[x index=0, y index=6] {\tableABS};
        \addplot table[x index=0, y index=7] {\tableABS};
        \addplot table[x index=0, y index=8] {\tableABS};
        
        \addplot[rngLine] coordinates {(1e0, 1e-20) (1e0, 1e0)};
        \addplot[rngLine] coordinates {(1e2, 1e-20) (1e2, 1e0)};
      \end{loglogaxis}
    \end{tikzpicture}};
    
  \node(rel)[right = 0cm of abs.south east, anchor = south west]{
    \begin{tikzpicture}
      \begin{loglogaxis}[%
        width  = .189\textwidth,
        height = .16\textheight,
        scale only axis,
        xmin = 1e-4,
        xmax = 1e+4,
        ymin = 1e-10,
        ymax = 1e+0,
        xminorticks = false,
        xlabel={\small Frequency (Hz)},
        xlabel style = {yshift = .3em},
        ylabel = {\small $\lVert H(j\omega) - \hH(j\omega) \rVert_{2} %
          / \lVert H(j\omega) \rVert_{2}$},
        ylabel style = {yshift = -.5em}]
        
        \pgfplotsset{cycle list shift=1}
        \addplot table[x index=0, y index=1] {\tableREL};
        \addplot table[x index=0, y index=2] {\tableREL};
        \addplot table[x index=0, y index=3] {\tableREL};
        \addplot table[x index=0, y index=4] {\tableREL};
        \addplot table[x index=0, y index=5] {\tableREL};
        \addplot table[x index=0, y index=6] {\tableREL};
        \addplot table[x index=0, y index=7] {\tableREL};
        \addplot table[x index=0, y index=8] {\tableREL};
        
        \addplot[rngLine] coordinates {(1e0, 1e-10) (1e0, 1e0)};
        \addplot[rngLine] coordinates {(1e2, 1e-10) (1e2, 1e0)};
      \end{loglogaxis}
    \end{tikzpicture}};
    
  \node[below = 0cm of tf.south west, anchor = north west]
    {(a) Sigma plot.};
  \node[below = 0cm of abs.south west, anchor = north west]
    {(b) Absolute errors.};
  \node[below = 0cm of rel.south west, anchor = north west]
    {(c) Relative errors.};
    
  \draw[draw = none](tf.west) -- (rel.east) node[midway](tmp){};
  \node(leg)[below = 3.3cm of tmp.center]{
    \begin{tikzpicture}
      \begin{axis}[%
        hide axis,
        scale only axis,
        width = 1mm,
        legend columns = 5, 
        legend style = {
          at     = {(0,0)},
          anchor = center,
          /tikz/every even column/.append style = {column sep = 0.5cm}}]
        \pgfplotsinvokeforeach{1,...,9}{\addplot coordinates {(0,0)};}
        
        \addlegendentry{Original};
        \addlegendentry{\myP};
        \addlegendentry{\myPM};
        \addlegendentry{\myPV};
        \addlegendentry{\myVP};
        \addlegendentry{\myVPM};
        \addlegendentry{\myV};
        \addlegendentry{\myFV};
        \addlegendentry{\mySO};
      \end{axis}
    \end{tikzpicture}};
\end{tikzpicture}

%% file: graphics/singlechain_sotlbt_a0.tikz
\begin{tikzpicture}[every axis/.append style={cycle list name = timelist}]
  \pgfplotstableread{graphics/data/singlechain_sotlbt_a0_step_abserr.dat}%
    \tableABSSTEP
  \pgfplotstableread{graphics/data/singlechain_sotlbt_a0_step_relerr.dat}%
    \tableRELSTEP
  \pgfplotstableread{graphics/data/singlechain_sotlbt_a0_sin_abserr.dat}%
    \tableABSSIN
  \pgfplotstableread{graphics/data/singlechain_sotlbt_a0_sin_relerr.dat}%
    \tableRELSIN
  
  \node(abs1){
    \begin{tikzpicture}
      \begin{semilogyaxis}[%
        width  = .32\textwidth,
        height = .15\textheight,
        scale only axis,
        xmin = 0,
        xmax = 100,
        ymin = 1e-10,
        ymax = 1e-2,
        xminorticks = false,
        yminorticks = false,
        xlabel={\small Time $t$ (sec)},
        xlabel style = {yshift = .3em},
        ylabel = {\small $\lVert y(t) - \hy(t) \rVert_{2}$},
        ylabel style = {yshift = -.5em}]
        
        \pgfplotsset{cycle list shift=1}
        \addplot table[x index=0, y index=1] {\tableABSSTEP};
        \addplot table[x index=0, y index=2] {\tableABSSTEP};
        \addplot table[x index=0, y index=3] {\tableABSSTEP};
        \addplot table[x index=0, y index=4] {\tableABSSTEP};
        \addplot table[x index=0, y index=5] {\tableABSSTEP};
        \addplot table[x index=0, y index=6] {\tableABSSTEP};
        \addplot table[x index=0, y index=7] {\tableABSSTEP};
        \addplot table[x index=0, y index=8] {\tableABSSTEP};
        
        \addplot[rngLine] coordinates {(20, 1e-10) (20, 1e-2)};
      \end{semilogyaxis}
    \end{tikzpicture}};
    
  \node(rel1)[right = 0.5cm of abs1.south east, anchor = south west]{
    \begin{tikzpicture}
      \begin{semilogyaxis}[%
        width  = .32\textwidth,
        height = .15\textheight,
        scale only axis,
        xmin = 0,
        xmax = 100,
        ymin = 1e-7,
        ymax = 1e-2,
        xminorticks = false,
        yminorticks = false,
        xlabel={\small Time $t$ (sec)},
        xlabel style = {yshift = .3em},
        ylabel = {\small $\lVert y(t) - \hy(t) \rVert_{2} %
          / \lVert y(t) \rVert_{2}$},
        ylabel style = {yshift = -.5em}]
        
        \pgfplotsset{cycle list shift=1}
        \addplot table[x index=0, y index=1] {\tableRELSTEP};
        \addplot table[x index=0, y index=2] {\tableRELSTEP};
        \addplot table[x index=0, y index=3] {\tableRELSTEP};
        \addplot table[x index=0, y index=4] {\tableRELSTEP};
        \addplot table[x index=0, y index=5] {\tableRELSTEP};
        \addplot table[x index=0, y index=6] {\tableRELSTEP};
        \addplot table[x index=0, y index=7] {\tableRELSTEP};
        \addplot table[x index=0, y index=8] {\tableRELSTEP};
        
        \addplot[rngLine] coordinates {(20, 1e-7) (20, 1e-2)};
      \end{semilogyaxis}
    \end{tikzpicture}};
    
  \node(abs2)[below = -.2cm of abs1.south]{
    \begin{tikzpicture}
      \begin{semilogyaxis}[%
        width  = .32\textwidth,
        height = .15\textheight,
        scale only axis,
        xmin = 0,
        xmax = 100,
        ymin = 1e-10,
        ymax = 1e-4,
        xminorticks = false,
        yminorticks = false,
        xlabel={\small Time $t$ (sec)},
        xlabel style = {yshift = .3em},
        ylabel = {\small $\lVert y(t) - \hy(t) \rVert_{2}$},
        ylabel style = {yshift = -.5em}]
        
        \pgfplotsset{cycle list shift=1}
        \addplot table[x index=0, y index=1] {\tableABSSIN};
        \addplot table[x index=0, y index=2] {\tableABSSIN};
        \addplot table[x index=0, y index=3] {\tableABSSIN};
        \addplot table[x index=0, y index=4] {\tableABSSIN};
        \addplot table[x index=0, y index=5] {\tableABSSIN};
        \addplot table[x index=0, y index=6] {\tableABSSIN};
        \addplot table[x index=0, y index=7] {\tableABSSIN};
        \addplot table[x index=0, y index=8] {\tableABSSIN};
        
        \addplot[rngLine] coordinates {(20, 1e-10) (20, 1e-4)};
      \end{semilogyaxis}
    \end{tikzpicture}};
    
  \node(rel2)[right = 0.5cm of abs2.south east, anchor = south west]{
    \begin{tikzpicture}
      \begin{semilogyaxis}[%
        width  = .32\textwidth,
        height = .15\textheight,
        scale only axis,
        xmin = 0,
        xmax = 100,
        ymin = 1e-8,
        ymax = 1e-1,
        xminorticks = false,
        yminorticks = false,
        xlabel={\small Time $t$ (sec)},
        xlabel style = {yshift = .3em},
        ylabel = {\small $\lVert y(t) - \hy(t) \rVert_{2} %
          / \lVert y(t) \rVert_{2}$},
        ylabel style = {yshift = -.5em}]
        
        \pgfplotsset{cycle list shift=1}
        \addplot table[x index=0, y index=1] {\tableRELSIN};
        \addplot table[x index=0, y index=2] {\tableRELSIN};
        \addplot table[x index=0, y index=3] {\tableRELSIN};
        \addplot table[x index=0, y index=4] {\tableRELSIN};
        \addplot table[x index=0, y index=5] {\tableRELSIN};
        \addplot table[x index=0, y index=6] {\tableRELSIN};
        \addplot table[x index=0, y index=7] {\tableRELSIN};
        \addplot table[x index=0, y index=8] {\tableRELSIN};
        
        \addplot[rngLine] coordinates {(20, 1e-8) (20, 1e-1)};
      \end{semilogyaxis}
    \end{tikzpicture}};
    
  \node[right = 0cm of abs1.south east, anchor = north west, yshift = 1cm]
    (a){(a)};
  \node[below = 0cm of rel1.south east, anchor = north west, yshift = 1cm]
    (b){(b)};
  \node[right = 0cm of abs2.south east, anchor = north west, yshift = 1cm]
    (c){(c)};
  \node[below = 0cm of rel2.south east, anchor = north west, yshift = 1cm]
    (d){(d)};
    
  \draw[draw = none](abs2.west) -- (abs2.west -| d.east) node[midway](tmp){};
  \node(leg)[below = 2.3cm of tmp.center]{
    \begin{tikzpicture}
      \begin{axis}[%
        hide axis,
        scale only axis,
        width = 1mm,
        legend columns = 4, 
        legend style = {
          at     = {(0,0)},
          anchor = center,
          /tikz/every even column/.append style = {column sep = 0.5cm}}]
        \pgfplotsset{cycle list shift=1}
        \pgfplotsinvokeforeach{1,...,8}{\addplot coordinates {(0,0)};}
        
        \addlegendentry{\myP};
        \addlegendentry{\myPM};
        \addlegendentry{\myPV};
        \addlegendentry{\myVP};
        \addlegendentry{\myVPM};
        \addlegendentry{\myV};
        \addlegendentry{\myFV};
        \addlegendentry{\mySO};
      \end{axis}
    \end{tikzpicture}};
\end{tikzpicture}

%% file: graphics/crankshaft_soflbt_a0_01.tikz
\begin{tikzpicture}[every axis/.append style={cycle list name = freqlist}]
  \pgfplotstableread{graphics/data/crankshaft_tf.dat}\tableFOM
  \pgfplotstableread{graphics/data/crankshaft_soflbt_a0_01_tf.dat}\tableROM
  \pgfplotstableread{graphics/data/crankshaft_soflbt_a0_01_abserr.dat}\tableABS
  \pgfplotstableread{graphics/data/crankshaft_soflbt_a0_01_relerr.dat}\tableREL
  
  \node(tf){
    \begin{tikzpicture}
      \begin{loglogaxis}[%
        width  = .189\textwidth,
        height = .16\textheight,
        scale only axis,
        xmin = 1e+0,
        xmax = 1e+6,
        ymin = 1e-10,
        ymax = 1e-3,
        xminorticks = false,
        xlabel={\small Frequency (Hz)},
        xlabel style = {yshift = .3em},
        ylabel = {\small $\lVert H(j\omega) \rVert_{2}$},
        ylabel style = {yshift = -.5em}]
        
        \addplot table[x index=0, y index=1] {\tableFOM};
        \addplot table[x index=0, y index=1] {\tableROM};
        \addplot table[x index=0, y index=2] {\tableROM};
        \addplot table[x index=0, y index=3] {\tableROM};
        \addplot table[x index=0, y index=4] {\tableROM};
        \addplot table[x index=0, y index=5] {\tableROM};
        \addplot table[x index=0, y index=6] {\tableROM};
        \addplot table[x index=0, y index=7] {\tableROM};
        \addplot table[x index=0, y index=8] {\tableROM};
        
        \addplot[rngLine] coordinates {(4e3, 1e-10) (4e3, 1e-3)};
        \addplot[rngLine] coordinates {(2e4, 1e-10) (2e4, 1e-3)};
      \end{loglogaxis}
    \end{tikzpicture}};
    
  \node(abs)[right = 0cm of tf.south east, anchor = south west]{
    \begin{tikzpicture}
      \begin{loglogaxis}[%
        width  = .189\textwidth,
        height = .16\textheight,
        scale only axis,
        xmin = 1e+0,
        xmax = 1e+6,
        ymin = 1e-12,
        ymax = 1e-2,
        xminorticks = false,
        xlabel={\small Frequency (Hz)},
        xlabel style = {yshift = .3em},
        ylabel = {\small $\lVert H(j\omega) - \hH(j\omega) \rVert_{2}$},
        ylabel style = {yshift = -.5em}]
        
        \pgfplotsset{cycle list shift=1}
        \addplot table[x index=0, y index=1] {\tableABS};
        \addplot table[x index=0, y index=2] {\tableABS};
        \addplot table[x index=0, y index=3] {\tableABS};
        \addplot table[x index=0, y index=4] {\tableABS};
        \addplot table[x index=0, y index=5] {\tableABS};
        \addplot table[x index=0, y index=6] {\tableABS};
        \addplot table[x index=0, y index=7] {\tableABS};
        \addplot table[x index=0, y index=8] {\tableABS};
        
        \addplot[rngLine] coordinates {(4e3, 1e-12) (4e3, 1e-2)};
        \addplot[rngLine] coordinates {(2e4, 1e-12) (2e4, 1e-2)};
      \end{loglogaxis}
    \end{tikzpicture}};
    
  \node(rel)[right = 0cm of abs.south east, anchor = south west]{
    \begin{tikzpicture}
      \begin{loglogaxis}[%
        width  = .189\textwidth,
        height = .16\textheight,
        scale only axis,
        xmin = 1e+0,
        xmax = 1e+6,
        ymin = 1e-5,
        ymax = 1e+2,
        xminorticks = false,
        xlabel={\small Frequency (Hz)},
        xlabel style = {yshift = .3em},
        ylabel = {\small $\lVert H(j\omega) - \hH(j\omega) \rVert_{2} %
          / \lVert H(j\omega) \rVert_{2}$},
        ylabel style = {yshift = -.5em}]
        
        \pgfplotsset{cycle list shift=1}
        \addplot table[x index=0, y index=1] {\tableREL};
        \addplot table[x index=0, y index=2] {\tableREL};
        \addplot table[x index=0, y index=3] {\tableREL};
        \addplot table[x index=0, y index=4] {\tableREL};
        \addplot table[x index=0, y index=5] {\tableREL};
        \addplot table[x index=0, y index=6] {\tableREL};
        \addplot table[x index=0, y index=7] {\tableREL};
        \addplot table[x index=0, y index=8] {\tableREL};
        
        \addplot[rngLine] coordinates {(4e3, 1e-5) (4e3, 1e+2)};
        \addplot[rngLine] coordinates {(2e4, 1e-5) (2e4, 1e+2)};
      \end{loglogaxis}
    \end{tikzpicture}};
    
  \node[below = 0cm of tf.south west, anchor = north west]
    {(a) Sigma plot.};
  \node[below = 0cm of abs.south west, anchor = north west]
    {(b) Absolute errors.};
  \node[below = 0cm of rel.south west, anchor = north west]
    {(c) Relative errors.};
    
  \draw[draw = none](tf.west) -- (rel.east) node[midway](tmp){};
  \node(leg)[below = 3.3cm of tmp.center]{
    \begin{tikzpicture}
      \begin{axis}[%
        hide axis,
        scale only axis,
        width = 1mm,
        legend columns = 5, 
        legend style = {
          at     = {(0,0)},
          anchor = center,
          /tikz/every even column/.append style = {column sep = 0.5cm}}]
        \pgfplotsinvokeforeach{1,...,9}{\addplot coordinates {(0,0)};}
        
        \addlegendentry{Original};
        \addlegendentry{\myP};
        \addlegendentry{\myPM};
        \addlegendentry{\myPV};
        \addlegendentry{\myVP};
        \addlegendentry{\myVPM};
        \addlegendentry{\myV};
        \addlegendentry{\myFV};
        \addlegendentry{\mySO};
      \end{axis}
    \end{tikzpicture}};
\end{tikzpicture}

%% file: graphics/crankshaft_soflbt_hybrid_local_a0_01.tikz
\begin{tikzpicture}[every axis/.append style={cycle list name = freqlist}]
  \pgfplotstableread{graphics/data/crankshaft_tf.dat}\tableFOM
  \pgfplotstableread{graphics/data/crankshaft_soflbt_hybrid_local_a0_01_tf.dat}\tableROM
  \pgfplotstableread{graphics/data/crankshaft_soflbt_hybrid_local_a0_01_abserr.dat}\tableABS
  \pgfplotstableread{graphics/data/crankshaft_soflbt_hybrid_local_a0_01_relerr.dat}\tableREL
  
  \node(tf){
    \begin{tikzpicture}
      \begin{loglogaxis}[%
        width  = .189\textwidth,
        height = .16\textheight,
        scale only axis,
        xmin = 1e+0,
        xmax = 1e+6,
        ymin = 1e-10,
        ymax = 1e-3,
        xminorticks = false,
        xlabel={\small Frequency (Hz)},
        xlabel style = {yshift = .3em},
        ylabel = {\small $\lVert H(j\omega) \rVert_{2}$},
        ylabel style = {yshift = -.5em}]
        
        \addplot table[x index=0, y index=1] {\tableFOM};
        \addplot table[x index=0, y index=1] {\tableROM};
        \addplot table[x index=0, y index=2] {\tableROM};
        \addplot table[x index=0, y index=3] {\tableROM};
        \addplot table[x index=0, y index=4] {\tableROM};
        \addplot table[x index=0, y index=5] {\tableROM};
        \addplot table[x index=0, y index=6] {\tableROM};
        \addplot table[x index=0, y index=7] {\tableROM};
        \addplot table[x index=0, y index=8] {\tableROM};
        
        \addplot[rngLine] coordinates {(4e3, 1e-10) (4e3, 1e-3)};
        \addplot[rngLine] coordinates {(2e4, 1e-10) (2e4, 1e-3)};
      \end{loglogaxis}
    \end{tikzpicture}};
    
  \node(abs)[right = 0cm of tf.south east, anchor = south west]{
    \begin{tikzpicture}
      \begin{loglogaxis}[%
        width  = .189\textwidth,
        height = .16\textheight,
        scale only axis,
        xmin = 1e+0,
        xmax = 1e+6,
        ymin = 1e-12,
        ymax = 1e-2,
        xminorticks = false,
        xlabel={\small Frequency (Hz)},
        xlabel style = {yshift = .3em},
        ylabel = {\small $\lVert H(j\omega) - \hH(j\omega) \rVert_{2}$},
        ylabel style = {yshift = -.5em}]
        
        \pgfplotsset{cycle list shift=1}
        \addplot table[x index=0, y index=1] {\tableABS};
        \addplot table[x index=0, y index=2] {\tableABS};
        \addplot table[x index=0, y index=3] {\tableABS};
        \addplot table[x index=0, y index=4] {\tableABS};
        \addplot table[x index=0, y index=5] {\tableABS};
        \addplot table[x index=0, y index=6] {\tableABS};
        \addplot table[x index=0, y index=7] {\tableABS};
        \addplot table[x index=0, y index=8] {\tableABS};
        
        \addplot[rngLine] coordinates {(4e3, 1e-12) (4e3, 1e-2)};
        \addplot[rngLine] coordinates {(2e4, 1e-12) (2e4, 1e-2)};
      \end{loglogaxis}
    \end{tikzpicture}};
    
  \node(rel)[right = 0cm of abs.south east, anchor = south west]{
    \begin{tikzpicture}
      \begin{loglogaxis}[%
        width  = .189\textwidth,
        height = .16\textheight,
        scale only axis,
        xmin = 1e+0,
        xmax = 1e+6,
        ymin = 1e-5,
        ymax = 1e+2,
        xminorticks = false,
        xlabel={\small Frequency (Hz)},
        xlabel style = {yshift = .3em},
        ylabel = {\small $\lVert H(j\omega) - \hH(j\omega) \rVert_{2} %
          / \lVert H(j\omega) \rVert_{2}$},
        ylabel style = {yshift = -.5em}]
        
        \pgfplotsset{cycle list shift=1}
        \addplot table[x index=0, y index=1] {\tableREL};
        \addplot table[x index=0, y index=2] {\tableREL};
        \addplot table[x index=0, y index=3] {\tableREL};
        \addplot table[x index=0, y index=4] {\tableREL};
        \addplot table[x index=0, y index=5] {\tableREL};
        \addplot table[x index=0, y index=6] {\tableREL};
        \addplot table[x index=0, y index=7] {\tableREL};
        \addplot table[x index=0, y index=8] {\tableREL};
        
        \addplot[rngLine] coordinates {(4e3, 1e-5) (4e3, 1e+2)};
        \addplot[rngLine] coordinates {(2e4, 1e-5) (2e4, 1e+2)};
      \end{loglogaxis}
    \end{tikzpicture}};
    
  \node[below = 0cm of tf.south west, anchor = north west]
    {(a) Sigma plot.};
  \node[below = 0cm of abs.south west, anchor = north west]
    {(b) Absolute errors.};
  \node[below = 0cm of rel.south west, anchor = north west]
    {(c) Relative errors.};
    
  \draw[draw = none](tf.west) -- (rel.east) node[midway](tmp){};
  \node(leg)[below = 3.3cm of tmp.center]{
    \begin{tikzpicture}
      \begin{axis}[%
        hide axis,
        scale only axis,
        width = 1mm,
        legend columns = 5, 
        legend style = {
          at     = {(0,0)},
          anchor = center,
          /tikz/every even column/.append style = {column sep = 0.5cm}}]
        \pgfplotsinvokeforeach{1,...,9}{\addplot coordinates {(0,0)};}
        
        \addlegendentry{Original};
        \addlegendentry{\myP};
        \addlegendentry{\myPM};
        \addlegendentry{\myPV};
        \addlegendentry{\myVP};
        \addlegendentry{\myVPM};
        \addlegendentry{\myV};
        \addlegendentry{\myFV};
        \addlegendentry{\mySO};
      \end{axis}
    \end{tikzpicture}};
\end{tikzpicture}

%% file: graphics/crankshaft_sotlbt_mixed_a0_01.tikz
\begin{tikzpicture}[every axis/.append style={scaled x ticks = false,
  x tick label style = {/pgf/number format/fixed},
  cycle list name = timelist}]
  \pgfplotstableread{graphics/data/crankshaft_sotlbt_mixed_a0_01_step_abserr.dat}%
    \tableABSSTEP
  \pgfplotstableread{graphics/data/crankshaft_sotlbt_mixed_a0_01_step_relerr.dat}%
    \tableRELSTEP
  \pgfplotstableread{graphics/data/crankshaft_sotlbt_mixed_a0_01_sin_abserr.dat}%
    \tableABSSIN
  \pgfplotstableread{graphics/data/crankshaft_sotlbt_mixed_a0_01_sin_relerr.dat}%
    \tableRELSIN
  
  \node(abs1){
    \begin{tikzpicture}
      \begin{semilogyaxis}[%
        width  = .32\textwidth,
        height = .15\textheight,
        scale only axis,
        xmin = 0,
        xmax = 0.05,
        ymin = 1e-7,
        ymax = 1e-3,
        xminorticks = false,
        yminorticks = false,
        xlabel={\small Time $t$ (sec)},
        xlabel style = {yshift = .3em},
        ylabel = {\small $\lVert y(t) - \hy(t) \rVert_{2}$},
        ylabel style = {yshift = -.5em}]
        
        \pgfplotsset{cycle list shift=1}
        \addplot table[x index=0, y index=1] {\tableABSSTEP};
        \addplot table[x index=0, y index=2] {\tableABSSTEP};
        \addplot table[x index=0, y index=3] {\tableABSSTEP};
        \addplot table[x index=0, y index=4] {\tableABSSTEP};
        \addplot table[x index=0, y index=5] {\tableABSSTEP};
        \addplot table[x index=0, y index=6] {\tableABSSTEP};
        \addplot table[x index=0, y index=7] {\tableABSSTEP};
        \addplot table[x index=0, y index=8] {\tableABSSTEP};
        
        \addplot[rngLine] coordinates {(0.01, 1e-7) (0.01, 1e-3)};
      \end{semilogyaxis}
    \end{tikzpicture}};
    
  \node(rel1)[right = 0.5cm of abs1.south east, anchor = south west]{
    \begin{tikzpicture}
      \begin{semilogyaxis}[%
        width  = .32\textwidth,
        height = .15\textheight,
        scale only axis,
        xmin = 0,
        xmax = 0.05,
        ymin = 1e-5,
        ymax = 1e-1,
        xminorticks = false,
        yminorticks = false,
        xlabel={\small Time $t$ (sec)},
        xlabel style = {yshift = .3em},
        ylabel = {\small $\lVert y(t) - \hy(t) \rVert_{2} %
          / \lVert y(t) \rVert_{2}$},
        ylabel style = {yshift = -.5em}]
        
        \pgfplotsset{cycle list shift=1}
        \addplot table[x index=0, y index=1] {\tableRELSTEP};
        \addplot table[x index=0, y index=2] {\tableRELSTEP};
        \addplot table[x index=0, y index=3] {\tableRELSTEP};
        \addplot table[x index=0, y index=4] {\tableRELSTEP};
        \addplot table[x index=0, y index=5] {\tableRELSTEP};
        \addplot table[x index=0, y index=6] {\tableRELSTEP};
        \addplot table[x index=0, y index=7] {\tableRELSTEP};
        \addplot table[x index=0, y index=8] {\tableRELSTEP};
        
        \addplot[rngLine] coordinates {(0.01, 1e-5) (0.01, 1e-1)};
      \end{semilogyaxis}
    \end{tikzpicture}};
    
  \node(abs2)[below = -.2cm of abs1.south]{
    \begin{tikzpicture}
      \begin{semilogyaxis}[%
        width  = .32\textwidth,
        height = .15\textheight,
        scale only axis,
        xmin = 0,
        xmax = 0.05,
        ymin = 1e-7,
        ymax = 1e-3,
        xminorticks = false,
        yminorticks = false,
        xlabel={\small Time $t$ (sec)},
        xlabel style = {yshift = .3em},
        ylabel = {\small $\lVert y(t) - \hy(t) \rVert_{2}$},
        ylabel style = {yshift = -.5em}]
        
        \pgfplotsset{cycle list shift=1}
        \addplot table[x index=0, y index=1] {\tableABSSIN};
        \addplot table[x index=0, y index=2] {\tableABSSIN};
        \addplot table[x index=0, y index=3] {\tableABSSIN};
        \addplot table[x index=0, y index=4] {\tableABSSIN};
        \addplot table[x index=0, y index=5] {\tableABSSIN};
        \addplot table[x index=0, y index=6] {\tableABSSIN};
        \addplot table[x index=0, y index=7] {\tableABSSIN};
        \addplot table[x index=0, y index=8] {\tableABSSIN};
        
        \addplot[rngLine] coordinates {(0.01, 1e-7) (0.01, 1e-3)};
      \end{semilogyaxis}
    \end{tikzpicture}};
    
  \node(rel2)[right = 0.5cm of abs2.south east, anchor = south west]{
    \begin{tikzpicture}
      \begin{semilogyaxis}[%
        width  = .32\textwidth,
        height = .15\textheight,
        scale only axis,
        xmin = 0,
        xmax = 0.05,
        ymin = 1e-5,
        ymax = 1e-1,
        xminorticks = false,
        yminorticks = false,
        xlabel={\small Time $t$ (sec)},
        xlabel style = {yshift = .3em},
        ylabel = {\small $\lVert y(t) - \hy(t) \rVert_{2} %
          / \lVert y(t) \rVert_{2}$},
        ylabel style = {yshift = -.5em}]
        
        \pgfplotsset{cycle list shift=1}
        \addplot table[x index=0, y index=1] {\tableRELSIN};
        \addplot table[x index=0, y index=2] {\tableRELSIN};
        \addplot table[x index=0, y index=3] {\tableRELSIN};
        \addplot table[x index=0, y index=4] {\tableRELSIN};
        \addplot table[x index=0, y index=5] {\tableRELSIN};
        \addplot table[x index=0, y index=6] {\tableRELSIN};
        \addplot table[x index=0, y index=7] {\tableRELSIN};
        \addplot table[x index=0, y index=8] {\tableRELSIN};
        
        \addplot[rngLine] coordinates {(0.01, 1e-5) (0.01, 1e-1)};
      \end{semilogyaxis}
    \end{tikzpicture}};
    
  \node[right = 0cm of abs1.south east, anchor = north west, yshift = 1cm]
    (a){(a)};
  \node[below = 0cm of rel1.south east, anchor = north west, yshift = 1cm]
    (b){(b)};
  \node[right = 0cm of abs2.south east, anchor = north west, yshift = 1cm]
    (c){(c)};
  \node[below = 0cm of rel2.south east, anchor = north west, yshift = 1cm]
    (d){(d)};
    
  \draw[draw = none](abs2.west) -- (abs2.west -| d.east) node[midway](tmp){};
  \node(leg)[below = 2.3cm of tmp.center]{
    \begin{tikzpicture}
      \begin{axis}[%
        hide axis,
        scale only axis,
        width = 1mm,
        legend columns = 4, 
        legend style = {
          at     = {(0,0)},
          anchor = center,
          /tikz/every even column/.append style = {column sep = 0.5cm}}]
        \pgfplotsset{cycle list shift=1}
        \pgfplotsinvokeforeach{1,...,8}{\addplot coordinates {(0,0)};}
        
        \addlegendentry{\myP};
        \addlegendentry{\myPM};
        \addlegendentry{\myPV};
        \addlegendentry{\myVP};
        \addlegendentry{\myVPM};
        \addlegendentry{\myV};
        \addlegendentry{\myFV};
        \addlegendentry{\mySO};
      \end{axis}
    \end{tikzpicture}};
\end{tikzpicture}

%% file: graphics/fish_tail_soflbt_a0.tikz
\begin{tikzpicture}[every axis/.append style={cycle list name = freqlist}]
  \pgfplotstableread{graphics/data/fish_tail_tf.dat}\tableFOM
  \pgfplotstableread{graphics/data/fish_tail_soflbt_a0_tf.dat}\tableROM
  \pgfplotstableread{graphics/data/fish_tail_soflbt_a0_abserr.dat}\tableABS
  \pgfplotstableread{graphics/data/fish_tail_soflbt_a0_relerr.dat}\tableREL
  
  \node(tf){
    \begin{tikzpicture}
      \begin{loglogaxis}[%
        width  = .189\textwidth,
        height = .16\textheight,
        scale only axis,
        xmin = 1e-4,
        xmax = 1e+4,
        ymin = 1e-12,
        ymax = 1e-2,
        xminorticks = false,
        xlabel={\small Frequency (Hz)},
        xlabel style = {yshift = .3em},
        ylabel = {\small $\lVert H(j\omega) \rVert_{2}$},
        ylabel style = {yshift = -.5em}]
        
        \addplot table[x index=0, y index=1] {\tableFOM};
        \addplot table[x index=0, y index=1] {\tableROM};
        \addplot table[x index=0, y index=2] {\tableROM};
        \addplot table[x index=0, y index=3] {\tableROM};
        \addplot table[x index=0, y index=4] {\tableROM};
        \addplot table[x index=0, y index=5] {\tableROM};
        \addplot table[x index=0, y index=6] {\tableROM};
        \addplot table[x index=0, y index=7] {\tableROM};
        \addplot table[x index=0, y index=8] {\tableROM};
        
        \addplot[rngLine] coordinates {(2e1, 1e-12) (2e1, 1e-2)};
      \end{loglogaxis}
    \end{tikzpicture}};
    
  \node(abs)[right = 0cm of tf.south east, anchor = south west]{
    \begin{tikzpicture}
      \begin{loglogaxis}[%
        width  = .189\textwidth,
        height = .16\textheight,
        scale only axis,
        xmin = 1e-4,
        xmax = 1e+4,
        ymin = 1e-12,
        ymax = 1e-2,
        xminorticks = false,
        xlabel={\small Frequency (Hz)},
        xlabel style = {yshift = .3em},
        ylabel = {\small $\lVert H(j\omega) - \hH(j\omega) \rVert_{2}$},
        ylabel style = {yshift = -.5em}]
        
        \pgfplotsset{cycle list shift=1}
        \addplot table[x index=0, y index=1] {\tableABS};
        \addplot table[x index=0, y index=2] {\tableABS};
        \addplot table[x index=0, y index=3] {\tableABS};
        \addplot table[x index=0, y index=4] {\tableABS};
        \addplot table[x index=0, y index=5] {\tableABS};
        \addplot table[x index=0, y index=6] {\tableABS};
        \addplot table[x index=0, y index=7] {\tableABS};
        \addplot table[x index=0, y index=8] {\tableABS};
        
        \addplot[rngLine] coordinates {(2e1, 1e-12) (2e1, 1e-2)};
      \end{loglogaxis}
    \end{tikzpicture}};
    
  \node(rel)[right = 0cm of abs.south east, anchor = south west]{
    \begin{tikzpicture}
      \begin{loglogaxis}[%
        width  = .189\textwidth,
        height = .16\textheight,
        scale only axis,
        xmin = 1e-4,
        xmax = 1e+4,
        ymin = 1e-4,
        ymax = 1e+2,
        xminorticks = false,
        xlabel={\small Frequency (Hz)},
        xlabel style = {yshift = .3em},
        ylabel = {\small $\lVert H(j\omega) - \hH(j\omega) \rVert_{2} %
          / \lVert H(j\omega) \rVert_{2}$},
        ylabel style = {yshift = -.5em}]
        
        \pgfplotsset{cycle list shift=1}
        \addplot table[x index=0, y index=1] {\tableREL};
        \addplot table[x index=0, y index=2] {\tableREL};
        \addplot table[x index=0, y index=3] {\tableREL};
        \addplot table[x index=0, y index=4] {\tableREL};
        \addplot table[x index=0, y index=5] {\tableREL};
        \addplot table[x index=0, y index=6] {\tableREL};
        \addplot table[x index=0, y index=7] {\tableREL};
        \addplot table[x index=0, y index=8] {\tableREL};
        
        \addplot[rngLine] coordinates {(2e1, 1e-4) (2e1, 1e+2)};
      \end{loglogaxis}
    \end{tikzpicture}};
    
  \node[below = 0cm of tf.south west, anchor = north west]
    {(a) Sigma plot.};
  \node[below = 0cm of abs.south west, anchor = north west]
    {(b) Absolute errors.};
  \node[below = 0cm of rel.south west, anchor = north west]
    {(c) Relative errors.};
    
  \draw[draw = none](tf.west) -- (rel.east) node[midway](tmp){};
  \node(leg)[below = 3.3cm of tmp.center]{
    \begin{tikzpicture}
      \begin{axis}[%
        hide axis,
        scale only axis,
        width = 1mm,
        legend columns = 5, 
        legend style = {
          at     = {(0,0)},
          anchor = center,
          /tikz/every even column/.append style = {column sep = 0.5cm}}]
        \pgfplotsinvokeforeach{1,...,9}{\addplot coordinates {(0,0)};}
        
        \addlegendentry{Original};
        \addlegendentry{\myP};
        \addlegendentry{\myPM};
        \addlegendentry{\myPV};
        \addlegendentry{\myVP};
        \addlegendentry{\myVPM};
        \addlegendentry{\myV};
        \addlegendentry{\myFV};
        \addlegendentry{\mySO};
      \end{axis}
    \end{tikzpicture}};
\end{tikzpicture}

%% file: graphics/fish_tail_sotlbt_a0.tikz
\begin{tikzpicture}[every axis/.append style={cycle list name = timelist}]
  \pgfplotstableread{graphics/data/fish_tail_sotlbt_a0_step_abserr.dat}%
    \tableABSSTEP
  \pgfplotstableread{graphics/data/fish_tail_sotlbt_a0_step_relerr.dat}%
    \tableRELSTEP
  \pgfplotstableread{graphics/data/fish_tail_sotlbt_a0_sin_abserr.dat}%
    \tableABSSIN
  \pgfplotstableread{graphics/data/fish_tail_sotlbt_a0_sin_relerr.dat}%
    \tableRELSIN
  
  \node(abs1){
    \begin{tikzpicture}
      \begin{semilogyaxis}[%
        width  = .32\textwidth,
        height = .15\textheight,
        scale only axis,
        xmin = 0,
        xmax = 2,
        ymin = 1e-8,
        ymax = 1e-2,
        xminorticks = false,
        yminorticks = false,
        xlabel={\small Time $t$ (sec)},
        xlabel style = {yshift = .3em},
        ylabel = {\small $\lVert y(t) - \hy(t) \rVert_{2}$},
        ylabel style = {yshift = -.5em}]
        
        \pgfplotsset{cycle list shift=1}
        \addplot table[x index=0, y index=1] {\tableABSSTEP};
        \addplot table[x index=0, y index=2] {\tableABSSTEP};
        \addplot table[x index=0, y index=3] {\tableABSSTEP};
        \addplot table[x index=0, y index=4] {\tableABSSTEP};
        \addplot table[x index=0, y index=5] {\tableABSSTEP};
        \addplot table[x index=0, y index=6] {\tableABSSTEP};
        \addplot table[x index=0, y index=7] {\tableABSSTEP};
        \addplot table[x index=0, y index=8] {\tableABSSTEP};
        
        \addplot[rngLine] coordinates {(0.5, 1e-8) (0.5, 1e-2)};
      \end{semilogyaxis}
    \end{tikzpicture}};
    
  \node(rel1)[right = 0.5cm of abs1.south east, anchor = south west]{
    \begin{tikzpicture}
      \begin{semilogyaxis}[%
        width  = .32\textwidth,
        height = .15\textheight,
        scale only axis,
        xmin = 0,
        xmax = 2,
        ymin = 1e-6,
        ymax = 1e+1,
        xminorticks = false,
        yminorticks = false,
        xlabel={\small Time $t$ (sec)},
        xlabel style = {yshift = .3em},
        ylabel = {\small $\lVert y(t) - \hy(t) \rVert_{2} %
          / \lVert y(t) \rVert_{2}$},
        ylabel style = {yshift = -.5em}]
        
        \pgfplotsset{cycle list shift=1}
        \addplot table[x index=0, y index=1] {\tableRELSTEP};
        \addplot table[x index=0, y index=2] {\tableRELSTEP};
        \addplot table[x index=0, y index=3] {\tableRELSTEP};
        \addplot table[x index=0, y index=4] {\tableRELSTEP};
        \addplot table[x index=0, y index=5] {\tableRELSTEP};
        \addplot table[x index=0, y index=6] {\tableRELSTEP};
        \addplot table[x index=0, y index=7] {\tableRELSTEP};
        \addplot table[x index=0, y index=8] {\tableRELSTEP};
        
        \addplot[rngLine] coordinates {(0.5, 1e-6) (0.5, 1e+1)};
      \end{semilogyaxis}
    \end{tikzpicture}};
    
  \node(abs2)[below = -.2cm of abs1.south]{
    \begin{tikzpicture}
      \begin{semilogyaxis}[%
        width  = .32\textwidth,
        height = .15\textheight,
        scale only axis,
        xmin = 0,
        xmax = 2,
        ymin = 1e-12,
        ymax = 1e-2,
        xminorticks = false,
        yminorticks = false,
        xlabel={\small Time $t$ (sec)},
        xlabel style = {yshift = .3em},
        ylabel = {\small $\lVert y(t) - \hy(t) \rVert_{2}$},
        ylabel style = {yshift = -.5em}]
        
        \pgfplotsset{cycle list shift=1}
        \addplot table[x index=0, y index=1] {\tableABSSIN};
        \addplot table[x index=0, y index=2] {\tableABSSIN};
        \addplot table[x index=0, y index=3] {\tableABSSIN};
        \addplot table[x index=0, y index=4] {\tableABSSIN};
        \addplot table[x index=0, y index=5] {\tableABSSIN};
        \addplot table[x index=0, y index=6] {\tableABSSIN};
        \addplot table[x index=0, y index=7] {\tableABSSIN};
        \addplot table[x index=0, y index=8] {\tableABSSIN};
        
        \addplot[rngLine] coordinates {(0.5, 1e-12) (0.5, 1e-2)};
      \end{semilogyaxis}
    \end{tikzpicture}};
    
  \node(rel2)[right = 0.5cm of abs2.south east, anchor = south west]{
    \begin{tikzpicture}
      \begin{semilogyaxis}[%
        width  = .32\textwidth,
        height = .15\textheight,
        scale only axis,
        xmin = 0,
        xmax = 2,
        ymin = 1e-8,
        ymax = 1e+2,
        xminorticks = false,
        yminorticks = false,
        xlabel={\small Time $t$ (sec)},
        xlabel style = {yshift = .3em},
        ylabel = {\small $\lVert y(t) - \hy(t) \rVert_{2} %
          / \lVert y(t) \rVert_{2}$},
        ylabel style = {yshift = -.5em}]
        
        \pgfplotsset{cycle list shift=1}
        \addplot table[x index=0, y index=1] {\tableRELSIN};
        \addplot table[x index=0, y index=2] {\tableRELSIN};
        \addplot table[x index=0, y index=3] {\tableRELSIN};
        \addplot table[x index=0, y index=4] {\tableRELSIN};
        \addplot table[x index=0, y index=5] {\tableRELSIN};
        \addplot table[x index=0, y index=6] {\tableRELSIN};
        \addplot table[x index=0, y index=7] {\tableRELSIN};
        \addplot table[x index=0, y index=8] {\tableRELSIN};
        
        \addplot[rngLine] coordinates {(0.5, 1e-8) (0.5, 1e+2)};
      \end{semilogyaxis}
    \end{tikzpicture}};
    
  \node[right = 0cm of abs1.south east, anchor = north west, yshift = 1cm]
    (a){(a)};
  \node[below = 0cm of rel1.south east, anchor = north west, yshift = 1cm]
    (b){(b)};
  \node[right = 0cm of abs2.south east, anchor = north west, yshift = 1cm]
    (c){(c)};
  \node[below = 0cm of rel2.south east, anchor = north west, yshift = 1cm]
    (d){(d)};
    
  \draw[draw = none](abs2.west) -- (abs2.west -| d.east) node[midway](tmp){};
  \node(leg)[below = 2.3cm of tmp.center]{
    \begin{tikzpicture}
      \begin{axis}[%
        hide axis,
        scale only axis,
        width = 1mm,
        legend columns = 4, 
        legend style = {
          at     = {(0,0)},
          anchor = center,
          /tikz/every even column/.append style = {column sep = 0.5cm}}]
        \pgfplotsset{cycle list shift=1}
        \pgfplotsinvokeforeach{1,...,8}{\addplot coordinates {(0,0)};}
        
        \addlegendentry{\myP};
        \addlegendentry{\myPM};
        \addlegendentry{\myPV};
        \addlegendentry{\myVP};
        \addlegendentry{\myVPM};
        \addlegendentry{\myV};
        \addlegendentry{\myFV};
        \addlegendentry{\mySO};
      \end{axis}
    \end{tikzpicture}};
\end{tikzpicture}